\documentclass[12pt,twoside]{article}

\usepackage[utf8]{inputenc}
\usepackage{amsmath,amsfonts,amssymb}
\usepackage{authblk}

\usepackage{fine-paper}
\usepackage{graphicx}
\usepackage{siunitx}
\usepackage{caption}
\usepackage[labelformat=simple]{subcaption}

\usepackage{nicefrac}
\usepackage{tikz}

\newcommand{\bS}{\boldsymbol{S}}			
\newcommand{\bPsi}{\boldsymbol{\psi}}	
\newcommand{\bPhi}{\boldsymbol{\phi}}	
\newcommand{\bu}{\boldsymbol{u}}			
\newcommand{\bE}{\boldsymbol{E}}			
\newcommand{\bW}{\boldsymbol{W}}			
\newcommand{\II}{I\! I}					
\newcommand{\III}{I\! I\! I}				

\newcommand{\diff}{\mathrm d}
\newcommand{\e}{\,\mathrm e}
\newcommand{\tr}{\mathrm{tr}}
\newcommand{\res}{\boldsymbol{\mathcal{R}}}
\newcommand{\jac}{\boldsymbol{\mathcal{J}}}
\newcommand{\proj}{\boldsymbol{\mathcal{P}}}

\newcommand{\eg}{\mbox{e.g.}}
\newcommand{\ie}{\mbox{i.e.}}
\newcommand{\etal}{\textit{et~al.\ }}

\title{The Variational Multiscale Formulation for the Fully-Implicit Log-Morphology Equation as a Tensor-Based Blood Damage Model}
\author{S. Haßler$^*$}
\author{L. Pauli}
\author{M. Behr}
\affil{Chair for Computational Analysis of Technical Systems (CATS)\\
    Center for Simulation and Data Science (JARA-CSD)\\
    RWTH Aachen University\\
    52056 Aachen, Germany\\
    \null\\
    $^*$hassler@cats.rwth-aachen.de}
\date{}

\begin{document}

\maketitle

\begin{abstract}
We derive a variational multiscale (VMS) finite element formulation for a viscoelastic, tensor-based blood damage model. The tensor equation is numerically stabilized by a logarithmic shape tensor description that prevents unphysical, negative eigenvalues. The resulting VMS stabilization terms for this so-called log-morph equation are presented together with their special numerical treatment. Results for a 2D rotating stirrer test case obtained from log-morph simulations with both SUPG and VMS stabilization show significantly improved numerical behavior if compared with Galerkin/least squares (GLS) stabilized untransformed morphology simulation results. The newly proposed method is also successfully applied to a state-of-the-art centrifugal ventricular assist device (VAD), and clear advantages of the VMS stabilization compared to the SUPG stabilized formulation are presented.

\end{abstract}

\section{Introduction}

Computational analysis has become a main step in the development of blood-handling medical devices. Computational Fluid Dynamics (CFD), in particular, helps to reduce the number of expensive and time-consuming experiments during prototype construction. CFD is able to predict the hydraulic performance and the flow behavior within such devices sufficiently and can therefore give useful information towards improvements of their design. This information is very valuable, since the seventh INTERMACS annual report \cite{Kirklin2015a} showed that adverse event rates for hemolysis, strokes and renal dysfunction have increased for patients with a ventricular assist device during the last survey period.

The accurate numerical prediction of hemolysis remains a challenging task. An empirical power law model is widely used to estimate the produced plasma-free hemoglobin in medical devices \cite{Fraser2011a,Yu2017a}. In such a model, the amount of generated plasma-free hemoglobin is proportional to the shear stress and the duration over which the RBCs are exposed to that stress. The power law parameters are normally determined by a fitting to experimental data obtained by Couette shearing device experiments using human or animal blood samples \cite{Heuser80,Wurzinger86a,Zhang2011a}. This modeling approach, the so-called stress-based hemolysis model, always computes a scalar measure of the shear stress from the flow field, which assumes an instantaneous deformation of the red blood cells.

Since such a stress-based model is not taking the viscoelastic behavior of RBCs into account, Arora~\etal~\cite{Arora2004a} proposed a strain-based model that estimates the droplet-like deformation of red blood cells in blood flow. The so-called morphology equation is able to account for the relaxation, elongation and rotation of the droplets in the flow. Arora~\etal~\cite{Arora2004a} used the model in a Lagrangian frame, and Pauli~\etal~\cite{Pauli2012b} proposed its application in an Eulerian frame.
Alternative viscoelastic blood damage models are introduced by different authors. Chen and Sharp \cite{Chen2011a} used a cell threshold model that they calibrated to channel flow experiments to estimate fatal hemolysis. However, such a threshold model cannot be applied to sublethal hemolysis. Arwatz and Smits \cite{Arwatz2013a} derived a scalar viscoelastic blood damage model, yet the generalisation to three-dimensional flows is not clear. Other strain-based approaches were proposed by Ezzeldin~\etal~\cite{Ezzeldin2015} and Sohrabi and Liu~\cite{Sohrabi2017a}. Ezzeldin \etal~\cite{Ezzeldin2015} compute the deformation of a high-fidelity RBC model based on a membrane energy composed of four contributions. Sohrabi and Liu~\cite{Sohrabi2017a} use a spring-connected network model for the description of the RBC membrane and use a pore formation model to estimate the generated plasma-free hemoglobin. Both works (Refs.~\cite{Ezzeldin2015,Sohrabi2017a}) describe the red blood cell membrane accurately but rely on a Lagrangian description, which can be biased by the choice of the tracers. Furthermore, it is computationally too costly to apply these methods to large-scale simulations such as those required for whole blood pumps.

In this study, we will use the morphology equation to compute the deformation of RBCs. Since the ellipsoidal RBC shape is described by the square root of the eigenvalues of the shape tensor, it is important that these eigenvalues are all positive, \ie, the shape tensor has to be positive-definite in order to describe a physical state. This positive-definiteness can be violated during the simulation, which is also a well-known issue for the related Oldroyd-B model. There, a logarithmic transformation of the conformation tensor fulfills the positive-definiteness by construction, and hence, improves the numerical stability of the method \cite{Fattal2004a,Knechtges2015a}. We will make use of such a transformation for the morphology equation and stabilize the transformed equations with a variational multiscale formulation.

The variational multiscale (VMS) method was introduced by Hughes \cite{Hughes1995a} and Hughes and Stewart \cite{Hughes96a} as a formulation to derive stabilization terms for the finite element method as analytical corrections from unresolved fine-scale contributions to the governing equations. In addition to the first applications to advection-diffusion and to Helmholtz equations, the VMS concept was applied to various other fields such as, \eg, the turbulence modeling for the incompressible Navier-Stokes equations \cite{Bazilevs2007a} and also to tensor-based equations such as the Oldroyd-B model \cite{Kwack2017a}.

The structure of the paper is as follows: In the next section, we will shortly introduce the morphology equation for the simulation of RBC deformation. We will motivate and present the transformation to a logarithmic shape tensor for enhanced numerical stability and its weak form for a finite-element implementation. The section will be completed with a detailed derivation of the variational multiscale formulation for this log-morph equation. Section \ref{sec:NumericalImplementation} will cover the numerical implementation of the newly derived VMS terms as well as the derivation and treatment of their directional derivatives for the usage in a Newton-Raphson solver. Two test cases, a two dimensional stirrer and a state-of-the-art ventricular assist device, will be investigated in Section \ref{sec:Results} and the advantages of the proposed log-morph formulation with VMS stabilization terms will be presented.

\section{Morphology Equation}

The behavior of red blood cells (RBCs) in blood flow is dependent on the shear rates in the fluid. At low shear rates, RBCs tend to form stack-like structures, called rouleaux \cite{Chien70a,Merrill66a}. A moderate increase of the fluid shear to tens of $\si{s^{-1}}$ breaks these rouleaux and the individual, biconcave RBCs move and tumble in the plasma \cite{Qin98a}. A further increase of the shear rate lets the RBCs rotate in the flow and deform to an ellipsoidal shape with a strained membrane \cite{Schmid-Schoenbein69a}; at around $\SI{4000}{s^{-1}}$, pores form in the lipid bilayer to release the membrane stresses, through which hemoglobin is released to the blood plasma \cite{Sohrabi2017a,Vitale2014}. Very high and unphysiological shear rates of about $\SI{42000}{s^{-1}}$ can ultimately lead to fatal hemolysis, the complete rupture of the RBC \cite{Leverett72a}. A recent study showed that the rotation of the RBC's membrane around the enclosed cytoplasm, the so-called tank treading motion, may not be observed for physiological viscosity ratios \cite{Lanotte2016a}.

\subsection{Droplet Model}

Arora \etal \cite{Arora2004a} proposed in 2004 a droplet model to account for the relaxation, deformation, and rotation of the RBCs, the so-called morphology model. This tensor-based model is able to describe the deformed, ellipsoidal shape of the RBCs. Arora \etal\ used this model in a Lagrangian frame, and included the tank treading motion in their model. Recent advances (cf. Pauli \etal~\cite{Pauli2012b}) used the tensor-based model in an Eulerian frame without the consideration of the tank-treading motion.

In the morphology model, the ellipsoidal RBC is described by a symmetric, positive-definite $3\times3$ matrix $\bS$. The residual $\res$ of the governing equation for the behavior of the RBCs in an external flow field $\bu$ is given by
\begin{equation}
\res\!\left(\bS\right) = \frac{\partial\bS}{\partial t} + \left(\bu\cdot\nabla\right)\bS
+ \underbrace{\alpha_{1}\left(\bS-g\!\left(\bS\right)\boldsymbol{1}\right)}_{\textrm{relaxation}}
- \underbrace{\alpha_{2}\left(\bE\bS+\bS\bE\right)}_{\textrm{elongation}}
- \underbrace{\alpha_{3}\left(\bW\bS-\bS\bW\right)}_{\textrm{rotation}} = \boldsymbol{0}.\label{eq:Res-morph}
\end{equation}
with unit matrix $\boldsymbol{1}$, the strain rate tensor $\bE = \left(\nabla\bu + \nabla\bu^{T}\right)/2$ and the vorticity tensor $\bW = \left(\nabla\bu - \nabla\bu^{T}\right)/2$. The parameter $g\!\left(\bS\right)=\frac{3\III\left(\bS\right)}{\II\left(\bS\right)}$,
with the third invariant $\III\!\left(\bS\right) = \det\!\left(\bS\right)$ and the second invariant $\II\!\left(\bS\right) = \left(\tr\!\left(\bS\right)^2 - \tr\!\left(\bS^2\right)\right) / 2$,
ensures the conservation of the volume of the RBCs. The parameters $\alpha_1 = \SI{5}{s^{-1}}$ and $\alpha_2 = \alpha_3 = \num{4.2298e-4}$ were derived by Arora \etal \cite{Arora2004a} from RBC relaxation and deformation properties.

The square roots of the eigenvalues of $\bS$ are the semi-axes lengths of the ellipsoid. For a computed shape $\bS$, one can compute the distortion $D$, a measure of the RBC's deformation, with the longest and shortest semi-axis, $L$ and $B$ by
\begin{equation}
D = \frac{L-B}{L+B}.
\end{equation}
One can show for a simple shear flow that the distortion is a function of the scalar shear stress $\sigma_{\textrm{f}} = 2 \mu \sqrt{-\II\!\left(\bE\right)}$, where $\mu$ is the blood viscosity. By inverting this relation, we get a formula for an effective shear stress that is acting on the RBC,
\begin{equation}
\sigma_{\textrm{eff}} = \frac{2 \mu \alpha_1 D}{\left(1 - D^2\right) \alpha_2},
\end{equation}
which is generalized and also used for complex flow and transient situations.

\subsection{Fully-Implicit Log-Morphology Formulation}

The positive definiteness of the morphology tensor $\bS$ is a necessary condition to describe a physical ellipsoidal shape. Nevertheless, this condition can be numerically violated leading to a diverged simulation. The morphology equation resembles the upper-convected Maxwell model, for which the positive-definiteness of the conformation tensor is also hard to satisfy numerically. Fattal and Kupferman \cite{Fattal2004a} proposed to use a log-conformation tensor --- a transformation that satisfies the positive definiteness by design --- to circumvent this problem. Knechtges \etal \cite{Knechtges2014a} and Knechtges \cite{Knechtges2015a} derived a fully-implicit log-conformation formulation, which can be also used with some modifications for the morphology equation.

The governing equation for the log-morphology tensor $\bPsi = \log\left(\bS\right)$, or in short, log-morph equation, can be derived as
\begin{align}
\res\!\left(\bPsi\right) = \frac{\partial\bPsi}{\partial t} + \left(\bu\cdot\nabla\right)\bPsi
& + \alpha_{1}\left(\boldsymbol{1} - g\!\left(\bPsi\right)\e^{-\bPsi}\right)
  - \alpha_{2}\boldsymbol{F}\!\left(\bPsi,\bE\right)
  - \alpha_{3}\left(\bW\bPsi - \bPsi\bW\right) = \boldsymbol{0},\label{eq:Res-log-morph}\\
\textrm{with }\boldsymbol{F}\!\left(\bPsi,\bE\right) =
& \frac{1}{\left(2\pi i\right)^{2}} \int_{\Gamma}\int_{\Gamma} f\!\left(z-z'\right) \frac{1}{z-\bPsi} \bE \frac{1}{z'-\bPsi} \diff z \diff z', \nonumber
\end{align}
representing a double Cauchy integral with the function $f\!\left(z\right) = \frac{z}{\tanh\left(z/2\right)}$ and a suitable contour $\Gamma$, which only encloses the complete spectrum of $\bPsi$, but not the poles of $f\!\left(z-z'\right)$. The derivation of the volume conservation term $g\!\left(\bPsi\right) = \frac{3}{\tr\left(\exp\left(-\bPsi\right)\right)}$
is shown in Appendix~\ref{sec:VolumeConservation}.

The arguments of Knechtges \cite{Knechtges2015a} hold also for the log-morph equation \eqref{eq:Res-log-morph} and a solution of this equation, transformed back, is also a solution of the untransformed morphology equation \eqref{eq:Res-morph}.

\subsubsection*{Weak Form}

We solve the log-morph equation using a space-time finite element method. To this end, we extrude our spatial domain in the time direction and divide the simulation time $t \in \left[0,T\right)$ into subintervals $I_n = \left[t_{n-1},t_n\right)$. We introduce space-time slabs $Q_n$ that are bounded by the spatial domain $\Omega_{n-1}$ at the time-step $t_{n-1}$ and by $\Omega_n$ at $t_n$. The spatial boundary of the space-time slab is denoted by $P_n$. Let the trial solution and weighting function spaces over $Q_n$ be given as $\mathcal{S}_n$ and $\mathcal{V}_n$. The variational form of the log-morph equation on a space-time slab $Q_n$ can then be given as: Find $\bPsi \in \mathcal{S}_n$ for the given initial condition $\bPsi\!\left(t_0^-\right) = \bPsi_0$ such that $\forall \, \bPhi \in \mathcal{V}_n$
\begin{equation}
0 = \underbrace{\int_{Q_n} \bPhi : \res\!\left(\bPsi\right) \diff Q}_{\equiv W\left(\bPhi;\bPsi\right)} + \underbrace{\int_{\Omega_{n-1}} \bPhi\!\left(t_{n-1}^+\right) : \left( \bPsi\!\left(t_{n-1}^+\right) - \bPsi\!\left(t_{n-1}^-\right) \right) \diff\Omega}_{\equiv D\left(\bPhi;\bPsi\right)} \label{eq:weak-log-morph}
\end{equation}
is satisfied. The second integral $D\!\left(\bPhi;\bPsi\right)$ is a discontinuous Galerkin term that weakly imposes the continuity of $\bPsi$ across space-time slabs. The Galerkin term $W\!\left(\bPhi;\bPsi\right)$ can be decomposed into its linear and nonlinear parts
\begin{align}
W\!\left(\bPhi;\bPsi\right) = &\overbrace{\int_{Q_n} \bPhi : \left[ \frac{\partial \bPsi}{\partial t} + \left(\bu \cdot  \nabla\right) \bPsi + \alpha_1 \boldsymbol{1} - \alpha_3 \left(\bW\bPsi - \bPsi\bW\right)\right] \diff Q}^{W_l\left(\bPhi;\bPsi\right)}  \nonumber\\
& \underbrace{-\int_{Q_n} \bPhi : \left[ \alpha_1 \frac{3 \exp\!\left(-\bPsi\right)}{\tr\!\left(\exp\!\left(-\bPsi\right)\right)} + \alpha_2 \boldsymbol{F}\!\left(\bPsi,\bE\right) \right] \diff Q}_{W_n\left(\bPhi;\bPsi\right)}.
\end{align}

\subsection{Variational Multiscale Formulation for Log-Morph}

In the VMS approach, one assumes that the real solution $\boldsymbol{U}$ to a general problem is composed of a coarse-scale solution $\boldsymbol{U}^h$ that can be resolved by the numerical method and an unresolved, fine-scale solution $\widetilde{\boldsymbol{U}}$ (cf. References~\cite{Hughes1995a,Bazilevs2007a}), \ie,
\begin{equation}
\boldsymbol{U} = \boldsymbol{U}^h + \widetilde{\boldsymbol{U}}.
\end{equation}
A further assumption is that the trial solution and weighting function spaces can be described as a direct sum of the coarse scale and fine scale spaces (cf. Ref.~\cite{Bazilevs2007a}), \ie,
\begin{equation}
\mathcal{S}_n = \mathcal{S}^h_n \oplus \widetilde{\mathcal{S}}_n \quad \textrm{and}\quad
\mathcal{V}_n = \mathcal{V}^h_n \oplus \widetilde{\mathcal{V}}_n.
\end{equation}

For the log-morph equation, we define the corresponding coarse scale trial solution and weighting function spaces as
\begin{align}
\mathcal{S}^h_n &= \left\{ \bPsi^h \in \left(C^0\!\left(\overline{Q_n}\right)\right)^6 \,\middle|\, \bPsi^h|_{P_n^D} = \boldsymbol{g}_{\bPsi} \right\}, \\
\mathcal{V}^h_n &= \left\{ \bPhi^h \in \left(C^0\!\left(\overline{Q_n}\right)\right)^6 \,\middle|\, \bPhi^h|_{P_n^D} = \boldsymbol{0} \right\},
\end{align}
on the finite element mesh $\overline{Q_n}$ and with the subset $P_n^D$ of the spatial space-time slab boundary $P_n$ where Dirichlet boundary conditions $\bPsi = \boldsymbol{g}_{\bPsi}$ are prescribed. Let us further assume that the fine scale trial functions are $\widetilde{\bPsi} \in \widetilde{\mathcal{S}}_n$ and the corresponding weighting functions are $\widetilde{\bPhi} \in \widetilde{\mathcal{V}}_n$ and that the VMS assumption
\begin{equation}
\bPsi = \bPsi^h + \widetilde{\bPsi}
\end{equation}
holds. The weak form of the log-morph equation \eqref{eq:weak-log-morph} is then decomposed into a set of coupled coarse and fine scale equations
\begin{align}
0 &= W_l\!\left(\bPhi^h;\bPsi^h + \widetilde{\bPsi}\right) + W_n\!\left(\bPhi^h;\bPsi^h + \widetilde{\bPsi}\right) + D\!\left(\bPhi^h;\bPsi^h + \widetilde{\bPsi}\right), \label{eq:VMS-coarse}\\
0 &= W_l\!\left(\widetilde{\bPhi}\,;\bPsi^h + \widetilde{\bPsi}\right) + W_n\!\left(\widetilde{\bPhi}\,;\bPsi^h + \widetilde{\bPsi}\right) + D\!\left(\widetilde{\bPhi}\,;\bPsi^h + \widetilde{\bPsi}\right). \label{eq:VMS-fine}
\end{align}
We can linearize the nonlinear terms in eqs.~\eqref{eq:VMS-coarse} and \eqref{eq:VMS-fine} using Fréchet derivatives \cite{Al-Mohy2009}:
\begin{equation}
W_n\!\left(\bullet;\bPsi^h + \widetilde{\bPsi}\right) \approx W_n\!\left(\bullet;\bPsi^h\right) + \left.\frac{\partial}{\partial\epsilon} W_n\!\left(\bullet;\bPsi^h + \epsilon \widetilde{\bPsi}\right)\right|_{\epsilon=0}, \label{eq:Frechet-derivative}
\end{equation}
which lets us rearrange eq.~\eqref{eq:VMS-fine} to
\begin{equation}
W_l\!\left(\widetilde{\bPhi};\widetilde{\bPsi}\right) + \left.\frac{\partial}{\partial\epsilon} W_n\!\left(\widetilde{\bPhi};\bPsi^h + \epsilon \widetilde{\bPsi}\right)\right|_{\epsilon=0} + D\!\left(\widetilde{\bPhi};\widetilde{\bPsi}\right) = -W_l\!\left(\widetilde{\bPhi};\bPsi^h\right) - W_n\!\left(\widetilde{\bPhi};\bPsi^h\right) - D\!\left(\widetilde{\bPhi};\bPsi^h\right). \label{eq:VMS-fine-separated}
\end{equation}
This is the weak form of a differential equation for the fine-scale solution $\widetilde{\bPsi}$, where the left-hand-side only contains operators acting on $\widetilde{\bPsi}$ and where only terms containing $\bPsi^h$ and $\res^h$ occur on the right-hand-side. As in Bazilevs \etal \cite{Bazilevs2007a}, eq.~\eqref{eq:VMS-fine-separated} hence tells us that the fine scale solution is a functional of the coarse scale solution and its residual
\begin{equation}
\widetilde{\bPsi} = \boldsymbol{\mathcal{F}}\!\left(\bPsi^h,\res\!\left(\bPsi^h\right)\right).
\end{equation}
Using physical reasoning, one can use a perturbation series for the fine scale solution that is dependent on powers of the coarse scale residual. One can show for the lowest order approximation that
\begin{equation}
\widetilde{\bPsi} \approx -\tau \res\!\left(\bPsi^h\right) \equiv -\tau \res^h, \label{eq:fine-approximation}
\end{equation}
with the stabilization parameter $\tau$ that we choose according to Shakib \etal \cite{Shakib91b} as
\begin{equation}
\tau = \alpha_\tau \left(\left(\frac{2}{\Delta t}\right)^2 + \bu \cdot \boldsymbol{G} \bu + \left\|\boldsymbol{L}\right\|_2\right)^{-\frac{1}{2}}, \label{eq:tau}
\end{equation}
\enlargethispage{5mm}
where $\boldsymbol{G}_{ij} = \sum_k \frac{\partial \xi^k}{\partial x^i} \frac{\partial \xi^k}{\partial x^j}$ is the covariant metric tensor mapping to a symmetric\footnote{We include the mapping to an equilateral triangle or a regular tetrahedron from Pauli \cite{Pauli2016b} in the definition of the metric tensor.} reference element \cite{Pauli2016b} and $\left\|\boldsymbol{L}\right\|_2$ is the spectral norm of the linearized source term, which we choose as the Jacobian of the source term (see eq.~\eqref{eq:jacobian}). We also use a scale factor $\alpha_\tau$ for our numerical implementation. In general, Hughes \cite{Hughes1995a} and others \cite{Hughes96a,Hughes98a,Bazilevs2007a} showed that the matrix-valued stabilization parameter $\boldsymbol{\tau}$ can be computed with the fine scales Green's operator.

Inserting the approximation from eq.~\eqref{eq:fine-approximation} into eq.~\eqref{eq:VMS-coarse} closes the system and lets us solve for the coarse scale log-morph tensor $\bPsi^h$. It can be easily deduced by an integration by parts that the VMS formulation includes the SUPG stabilization terms, since
\begin{align}
W_l\!\left(\bPhi^h; \bPsi^h - \tau\res^h\right) = W_l\!\left(\bPhi^h;\bPsi^h\right) &+ \int_{Q_n} \tau \left[ \frac{\partial \bPhi^h}{\partial t} + \left(\bu\cdot\nabla\right) \bPhi^h \right] : \res^h \, \diff Q \nonumber\\
&+ \int_{Q_n} \tau \alpha_3 \bPhi^h : \left[ \bW\res^h - \res^h\bW \right] \diff Q.
\end{align}
As mentioned before, the arising nonlinear terms have to be linearized with respect to the coarse scale solution. In order to compute the Fréchet derivatives, we have to use the following
relation,
\begin{equation}
\frac{\partial}{\partial\epsilon} \frac{1}{z-\epsilon\bPsi^h} = \frac{1}{z-\epsilon\bPsi^h} \bPsi^h \frac{1}{z-\epsilon\bPsi^h},
\end{equation}
and we need to compute the derivative of the matrix exponential for the $\alpha_1$-term,
\begin{align}
\left.\frac{\partial}{\partial\epsilon} \exp\!\left(-\bPsi^h+\epsilon\tau\res^h\right) \right|_{\epsilon=0} & = \frac{1}{2\pi i} \left.\frac{\partial}{\partial\epsilon} \int_{\Gamma} \e^{-z} \frac{1}{z - \bPsi^h + \epsilon\tau\res^h} \diff z \, \right|_{\epsilon=0} \nonumber \\
& = -\tau \underbrace{\frac{1}{2\pi i} \int_{\Gamma} \e^{-z} \frac{1}{z-\bPsi^h} \res^h \frac{1}{z-\bPsi^h} \diff z}_{\equiv \boldsymbol{K}\left(\bPsi^h,\res^h\right)}.
\end{align}
The resulting linearized terms $\boldsymbol{L}_{\alpha_1}$ and $\boldsymbol{L}_{\alpha_2}$ are:
\begin{equation}
\left. \frac{\partial}{\partial\epsilon} \frac{\exp\!\left(-\bPsi^h + \epsilon\tau\res^h\right)}{\tr\!\left(\exp\!\left(-\bPsi^h + \epsilon\tau\res^h\right)\right)} \right|_{\epsilon=0} = -\tau \underbrace{\left(\frac{\boldsymbol{K}\!\left(\bPsi^h,\res^h\right)}{\tr\!\left(\exp\!\left(-\bPsi^h\right)\right)} - \frac{\exp\!\left(-\bPsi^h\right)}{\tr\!\left(\exp\!\left(-\bPsi^h\right)\right)^{2}} \tr\!\left(\boldsymbol{K}\!\left(\bPsi^h,\res^h\right)\right)\right)}_{\equiv \boldsymbol{L}_{\alpha_1}\left(\bPsi^h,\res^h\right)},\label{eq:L_f1}
\end{equation}
and
\begin{gather}
\left. \frac{\partial}{\partial\epsilon} \boldsymbol{F}\!\left(\bPsi^h - \epsilon\tau\res^h, \bE\right)\right|_{\epsilon=0} \nonumber \\
=  -\tau \underbrace{\frac{1}{\left(2\pi i\right)^{2}} \int_{\Gamma} \int_{\Gamma} f\!\left(z-z'\right) \left[\frac{1}{z-\bPsi^h} \res^h \frac{1}{z-\bPsi^h} \bE \frac{1}{z'-\bPsi^h} + \frac{1}{z-\bPsi^h} \bE \frac{1}{z'-\bPsi^h} \res^h \frac{1}{z'-\bPsi^h}\right] \diff z \diff z'}_{\equiv \boldsymbol{L}_{\alpha_2}\left(\bPsi^h,\res^h,\bE\right)}.
\end{gather}
It is important to note, that the additional VMS terms are all traceless which ensures the volume conservation (cf. Appendix~\ref{sec:VolumeConservation}).

The resulting VMS stabilized equations are
\begin{align}
0 &= \int_{Q_n} \bPhi^h : \res\!\left(\bPsi^h\right) \diff Q + \int_{\Omega_{n-1}} \bPhi^h\!\left(t_{n-1}^+\right) : \left( \bPsi^h\!\left(t_{n-1}^+\right) - \bPsi^h\!\left(t_{n-1}^-\right) \right) \diff\Omega \nonumber\\
&+ \int_{Q_n} \tau \left[ \frac{\partial \bPhi^h}{\partial t} + \left(\bu\cdot\nabla\right) \bPhi^h \right] : \res\!\left(\bPsi^h\right) \diff Q + \int_{Q_n} \tau \alpha_3 \bPhi^h : \left[ \bW\res\!\left(\bPsi^h\right) - \res\!\left(\bPsi^h\right)\bW \right] \diff Q \nonumber\\
&+ \int_{Q_n} \tau \bPhi^h : \left[ 3 \alpha_1 \boldsymbol{L}_{\alpha_1}\!\left(\bPsi^h,\res\!\left(\bPsi^h\right)\right) + \alpha_2 \boldsymbol{L}_{\alpha_2}\!\left(\bPsi^h,\res\!\left(\bPsi^h\right),\bE\right) \right] \diff Q. \label{eq:weak-log-morph-vms}
\end{align}

\section{Numerical Implementation} \label{sec:NumericalImplementation}

For the numerical implementation of the log-morph equation with the VMS stabilization, we need to perform some further steps. The discretized equations are going to be solved using a Newton-Raphson algorithm combined with a GMRES solver for the resulting linearized equation system. The Newton-Raphson algorithm requires a further directional derivative with respect to $\delta\bPsi^h$ (also a Fréchet derivative) for the linearization of the system. This derivative has to be computed for every term arising in $W_l\!\left(\bPhi^h;\bPsi^h-\tau\res^h\right)$ and $W_n\!\left(\bPhi^h;\bPsi^h-\tau\res^h\right)$.

In order to evaluate the Cauchy integrals that arise, we use projectors $\proj_{i}=\boldsymbol{e}_{i}\boldsymbol{e}_{i}^{T}$ that project onto the one-dimensional subspaces of the eigenvalues $\lambda_{i}$ of $\bPsi^h$ spanned by the corresponding eigenvectors $\boldsymbol{e}_{i}$, \ie,
\begin{equation}
\frac{1}{z-\bPsi^h}=\sum_{i=1}^{d}\frac{1}{z-\lambda_{i}}\proj_{i}.
\end{equation}
This leads to Cauchy integrals that can be evaluated using the residue theorem, which are similar to
\begin{equation}
\frac{1}{2\pi i} \int_{\Gamma} m\!\left(z\right) \frac{1}{z-\lambda} \diff z = m\!\left(\lambda\right),
\end{equation}
with an arbitrary function $m\!\left(z\right)$ without a singularity at $z = \lambda$.

\subsection{Numerical Evaluation of the VMS Weak Form}

The evaluation of the Cauchy integrals is done according to Knechtges~\cite{Knechtges2015a}, whose numerical treatment of the arising prefactors can also be used for our equations. The resulting $\alpha_1$-term in the residual can be computed with the matrix exponential
\begin{equation}
\exp\!\left(-\bPsi^h\right) = \frac{1}{2\pi i} \int_\Gamma \e^{-z} \frac{1}{z-\bPsi^h} \diff z = \sum_{i=1}^{d} \e^{-\lambda_i} \proj_i,
\end{equation}
from which
\begin{equation}
\tr\!\left(\exp\!\left(-\bPsi^h\right)\right) = \sum_{i=1}^d \e^{-\lambda_i}
\end{equation}
immediately follows. The $\alpha_2$-term can be evaluated with
\begin{equation}
\boldsymbol{F}\!\left(\bPsi^h,\bE\right) = \sum_{i,j=1}^d f\!\left(\lambda_i - \lambda_j\right) \proj_i \bE \proj_j.
\end{equation}
The function $f(x) = \frac{x}{\tanh(x/2)}$ is set to 2 for $x = 0$ and can be readily used for values of $x \neq 0$.

Similar considerations yield the new VMS terms in eq.~\eqref{eq:weak-log-morph-vms}, where we have to further evaluate
\begin{equation}
\boldsymbol{K}\!\left(\bPsi^h;\res^h\right) = -\sum_{i,j=1}^d \e^{-\nicefrac{\lambda_{i}}{2}} \e^{-\nicefrac{\lambda_{j}}{2}} \frac{\sinh\!\left(\nicefrac{\left(\lambda_{i}-\lambda_{j}\right)}{2}\right)}{\nicefrac{\left(\lambda_{i}-\lambda_{j}\right)}{2}} \proj_{i}\res^h\proj_{j} 
\end{equation}
and take traces of the form
\begin{equation}
\tr\!\left(\proj_i \boldsymbol{A} \proj_j\right) = \boldsymbol{e}_i^T \boldsymbol{A} \boldsymbol{e}_i,
\end{equation}
with an arbitrary $3\times3$ matrix $\boldsymbol{A}$ for the determination of the $\boldsymbol{L}_{\alpha_1}\!\left(\bPsi^h,\res^h\right)$ term. We find for the other nonlinear VMS term
\begin{equation}
\boldsymbol{L}_{\alpha_2}\!\left(\bPsi^h,\res^h,\bE\right) =\sum_{i,j,k=1}^d \left(\frac{f\!\left(\lambda_{i}-\lambda_{k}\right) - f\!\left(\lambda_{j}-\lambda_{k}\right)}{\lambda_{i}-\lambda_{j}}\right) \left[\proj_{i}\res^h\proj_{j}\bE\proj_{k}+\proj_{k}\bE\proj_{j}\res^h\proj_{i}\right].
\end{equation}

As discussed by Knechtges~\cite{Knechtges2015a}, the prefactor $\frac{\sinh\left(x/2\right)}{(x/2)}$ is used for values $x \neq 0$, and replaced by 1 for $x = 0$. The prefactor in the $\boldsymbol{L}_{\alpha_2}$ term can be approximated by a Taylor series
\begin{equation}
\frac{f(x) - f(y)}{x - y} = f'\!\left(\frac{x + y}{2}\right) + \frac{\left(x - y\right)^2}{24} f'''\!\left(\frac{x + y}{2}\right) + \mathcal{O}\!\left(\left(x-y\right)^4\right) \label{eq:TaylorSmallDenominator}
\end{equation}
in the vicinity of small denominators, $|x - y| < 10^{-2}$. The derivatives are approximated by their Taylor series up to fourth order for small arguments $\left|\nicefrac{\left(x + y\right)}{2}\right| < 10^{-1}$.

\subsection{Numerical Evaluation of the Directional Derivatives for the Newton-Raphson Algorithm}

As mentioned before, we have to take a further derivative in the direction of $\delta\bPsi^h$ for the linearization in the Newton-Raphson algorithm. For the residual $\res^h$ this gives us the Jacobian
\begin{align}
\jac^h \equiv \jac\!\left(\bPsi^h;\delta\bPsi^h\right)
= &\left.\frac{\partial}{\partial\epsilon} \res\!\left(\bPsi^h + \epsilon\, \delta\bPsi^h\right)\right|_{\epsilon=0}
= \frac{\partial\delta\bPsi^h}{\partial t} + \left(\bu\cdot\nabla\right)\delta\bPsi^h \nonumber\\
 &\underbrace{- \alpha_1 3 \boldsymbol{L}_{\alpha_1}\!\left(\bPsi^h,\delta\bPsi^h\right)
  - \alpha_2 \boldsymbol{L}_{\alpha_2}\!\left(\bPsi^h,\delta\bPsi^h,\bE\right)
  - \alpha_3\left(\bW\delta\bPsi^h - \delta\bPsi^h\bW\right)}_{\equiv \boldsymbol{L}\left(\bPsi^h,\delta\bPsi^h,\bE,\bW\right)}, \label{eq:jacobian}
\end{align}
where $\delta\bPsi^h$ has to be used instead of $\tau\res^h$ in $\boldsymbol{L}_{\alpha_1}$ and $\boldsymbol{L}_{\alpha_2}$. This Jacobian is used for the Galerkin part as well as for the linearized VMS terms of the weak form in eq.~\eqref{eq:weak-log-morph-vms}.

For the nonlinear $\alpha_1$-term, we have to compute the directional derivative of $\boldsymbol{K}\!\left(\bPsi^h,\res^h\right)$ in the direction of $\delta\bPsi^h$, which leads to
\begin{align}
&\left.\frac{\partial}{\partial\epsilon} \boldsymbol{K}\!\left(\bPsi^h + \epsilon\,\delta\bPsi^h,\res\left(\bPsi^h + \epsilon\,\delta\bPsi^h\right)\right)\right|_{\epsilon=0}
= \frac{1}{2\pi i} \int_\Gamma \e^{-z} \frac{1}{z - \bPsi^h} \jac^h \frac{1}{z - \bPsi^h} \,\diff z  \nonumber\\
& +  \frac{1}{2\pi i} \int_\Gamma \e^{-z} \left[\frac{1}{z - \bPsi^h} \delta\bPsi^h \frac{1}{z - \bPsi^h} \res^h \frac{1}{z - \bPsi^h} + \frac{1}{z - \bPsi^h} \res^h \frac{1}{z - \bPsi^h} \delta\bPsi^h \frac{1}{z - \bPsi^h}\right] \diff z \nonumber\\
= & -\sum_{i,j=1}^d \e^{-\nicefrac{\lambda_i}{2}} \e^{-\nicefrac{\lambda_j}{2}} \frac{\sinh\!\left(\nicefrac{\left(\lambda_i-\lambda_j\right)}{2}\right)}{\nicefrac{\left(\lambda_i-\lambda_j\right)}{2}} \proj_i \jac^h \proj_j \nonumber \\
& +\sum_{i,j,k=1}^d \left(\frac{\e^{-\lambda_i}}{\left(\lambda_i - \lambda_j\right) \left(\lambda_i - \lambda_k\right)} + \frac{\e^{-\lambda_j}}{\left(\lambda_j - \lambda_i\right) \left(\lambda_j - \lambda_k\right)} + \frac{\e^{-\lambda_k}}{\left(\lambda_k - \lambda_i\right) \left(\lambda_k - \lambda_j\right)}\right) \nonumber\\
&\qquad \left[\proj_i \res^h \proj_j \delta\bPsi^h \proj_k + \proj_k \delta\bPsi^h \proj_j \res^h \proj_i\right].
\end{align}
The prefactor of the second term containing the residuals is also discussed by Knechtges~\cite{Knechtges2015a} and Taylor series approximations are used for small denominators $|x-y| < 10^{-3}$ and for the arising derivatives for arguments $\left|\nicefrac{(x + y)}{2}\right| < 10^{-3}$. For the complete directional derivative of the $\boldsymbol{L}_{\alpha_1}$-term in eq.~\eqref{eq:L_f1}, we have to consider the chain rule.

Accordingly, the directional derivative of $\boldsymbol{L}_{\alpha_2}$ becomes
\begin{align}
&\left. \frac{\partial}{\partial \epsilon} \boldsymbol{L}_{\alpha_2}\!\left(\bPsi^h + \epsilon\,\delta\bPsi^h,\res\!\left(\bPsi^h + \epsilon\, \delta\bPsi^h\right),\bE\right) \right|_{\epsilon = 0} \nonumber\\
= & \sum_{i,j,k=1}^d \frac{f\!\left(\lambda_i - \lambda_k\right) - f\!\left(\lambda_j - \lambda_k\right)}{\lambda_i - \lambda_j} \left[\proj_i \jac^h \proj_j \bE \proj_k + \proj_k \bE \proj_j \jac^h \proj_i\right]\nonumber \\
& +\sum_{i,j,k,l=1}^d \left(\frac{f\!\left(\lambda_i - \lambda_l\right)}{\left(\lambda_i - \lambda_j\right) \left(\lambda_i - \lambda_k\right)} + \frac{f\!\left(\lambda_j - \lambda_l\right)}{\left(\lambda_j - \lambda_i\right) \left(\lambda_j - \lambda_k\right)} + \frac{f\!\left(\lambda_k - \lambda_l\right)}{\left(\lambda_k - \lambda_i\right) \left(\lambda_k - \lambda_j\right)}\right) \nonumber \\
&\qquad \left[\proj_i \delta\bPsi^h \proj_j \res^h \proj_k \bE \proj_l + \proj_i \res^h \proj_j \delta\bPsi^h \proj_k \bE \proj_l + \proj_l \bE \proj_k \delta\bPsi^h \proj_j \res^h \proj_i + \proj_l \bE \proj_k \res^h \proj_j \delta\bPsi^h \proj_i\right] \nonumber \\
& +\sum_{i,j,k,l=1}^d \frac{f\!\left(\lambda_i - \lambda_k\right) - f\!\left(\lambda_j - \lambda_k\right) + f\!\left(\lambda_i - \lambda_l\right) - f\!\left(\lambda_j - \lambda_l\right)}{\left(\lambda_i - \lambda_j\right) \left(\lambda_k - \lambda_l\right)} \nonumber \\
&\qquad \left[\proj_i \res^h \proj_j \bE \proj_k \delta\bPsi^h \proj_l + \proj_l \delta\bPsi^h \proj_k \bE \proj_j \res^h \proj_i\right]. \label{eq:DerivativeLa2}
\end{align}
The first prefactor was already discussed in the previous section, but the second and third prefactors have to be treated with special care. A complete discussion of these prefactors is done in Appendix~\ref{sec:Prefactors}.

\section{Results} \label{sec:Results}

In order to test the newly proposed log-morph equation with VMS stabilization, we investigate two test cases: A simple two-dimensional rotating stirrer in a square box, and a state-of-the-art ventricular assist device (VAD). We consider three different discretizations for our test cases: the untransformed morphology equation with GLS stabilization (morph-GLS, we refer to Ref.~\cite{Pauli2016} for a detailed derivation), the log-morph equation with only SUPG stabilization (log-morph-SUPG), and the log-morph equation with full VMS stabilization (log-morph-VMS). The computations for this section were performed on the supercomputer JURECA at Forschungszentrum Jülich \cite{jureca}.

\subsection{Stirrer Test Case}\label{sec:Stirrer}

\begin{figure}
\begin{minipage}{0.48\textwidth}
\centering
\input{pics/StirrerGeometry.tex}%
\captionof{figure}{Geometry of the 2D stirrer test case with the MRF interface (dashed circle) and a part of the computational mesh (adopted from Ref.~\cite{Pauli2014a}).}
\label{fig:StirrerGeometry}
\end{minipage}\hfill
\begin{minipage}{0.48\textwidth}
\centering
\includegraphics[height=8cm]{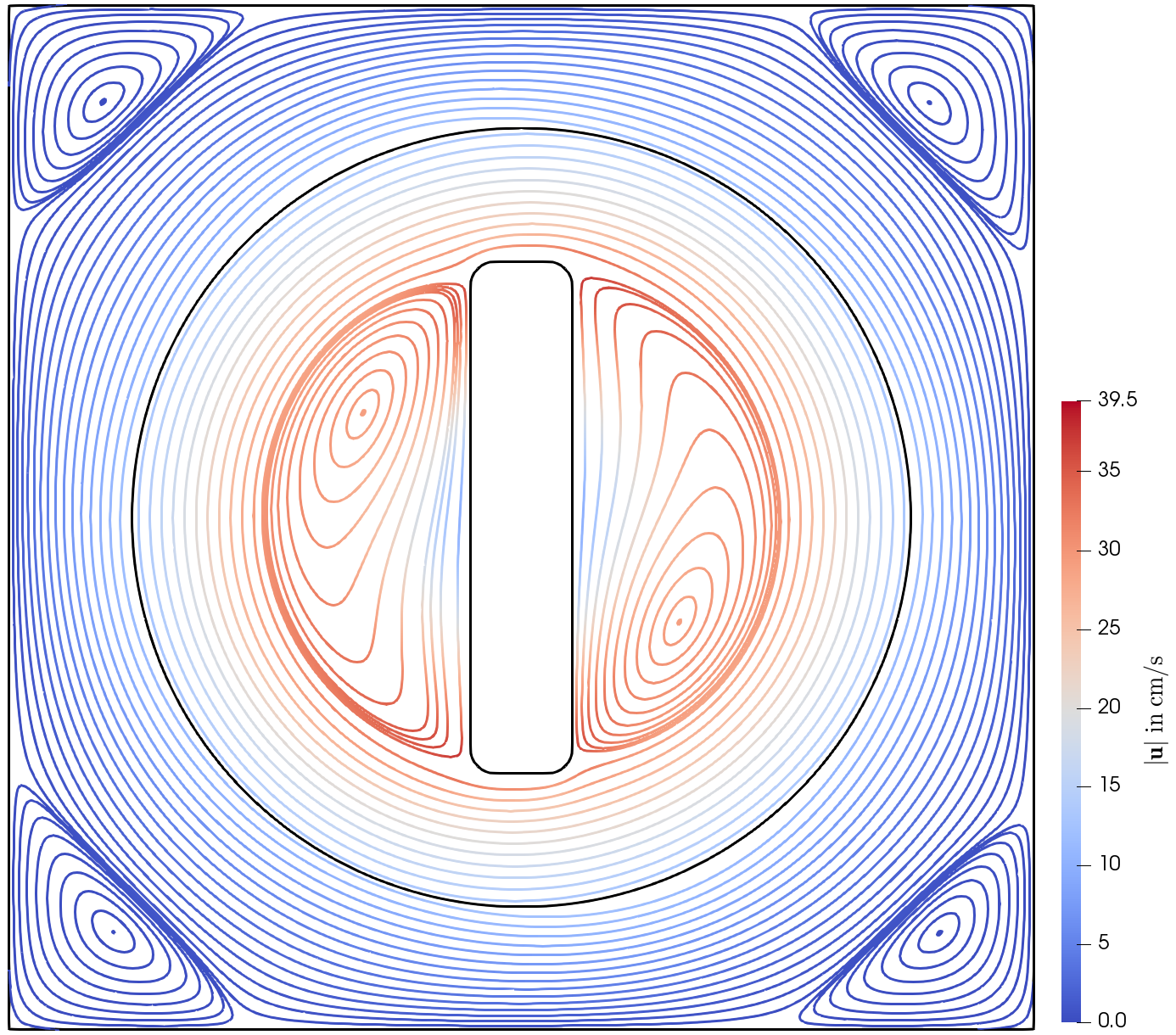}
\vspace{2mm}
\captionof{figure}{Velocity streamlines in the corresponding reference frame colored by the velocity magnitude in the inertial frame.\\\null}
\label{fig:StirrerVelocity}
\end{minipage}
\end{figure}

The computational mesh for the stirrer test case consists of $n_n = \num{46647}$ nodes and $n_e = \num{92262}$ unstructured triangular elements. A part of the mesh near the beam is shown in Fig. \ref{fig:StirrerGeometry} together with the dimensions of the geometry.
The flow solution is computed with the in-house deformable-spatial-domain/stabilized space-time finite element solver XNS. We use a Newtonian blood model with a density of $\rho = \SI{1054}{kg/m^3}$ and a viscosity of $\mu = \SI{0.0035}{Pa.s}$. We utilize the multiple reference frames (MRF) method (cf. Ref.~\cite{Pauli2014a}) to compute a steady-state approximation of the flow field for an angular velocity of $\omega = \SI{50\pi}{s^{-1}}$. The angular velocity is increased from 0 to the full rotational speed through the first 30 Newton-Raphson iterations $i_{\mathrm{NR}}$ according to
\begin{equation}
\omega = 50\pi \left(3 \, i_{\mathrm{NR}}^2 - 2 \, i_{\mathrm{NR}}^3\right).
\end{equation}
The steady-state flow field is converged after 34 iterations and is shown in Fig.~\ref{fig:StirrerVelocity}, where we show the streamlines in the corresponding reference frame.

The flow solution is used as an ambient velocity for the computation of the deformation of the RBCs. For this simple test case, we are able to choose a time step size $\Delta t = \SI{0.01}{s}$ for all three discretizations. We compute a quasi-steady deformation for an exposure time of $\SI{1}{s}$ using a semi-discrete time discretization scheme with a backward differentiation formula of second order (BDF2). Also for the morphology equation, we use an MRF approach. A Newton-Raphson algorithm combined with a preconditioned, restarted GMRES solver is used to solve the nonlinear system. For all three simulations we chose a Krylov space of $\num{10}$ and an ILUT-preconditioner with a maximal fill-in of $\num{20}$ and a threshold of $\num{e-4}$. In the morph-GLS case we have to further use an augmented Lagrangian method with a penalty parameter of $\varepsilon_p = \num{10000}$ in order to penalize deviations from the initial droplet volume (cf. Ref.~\cite{Pauli2016}).

\begin{figure}[tb]
\centering
\begin{subfigure}{.48\textwidth}
\centering
\includegraphics[width=\textwidth]{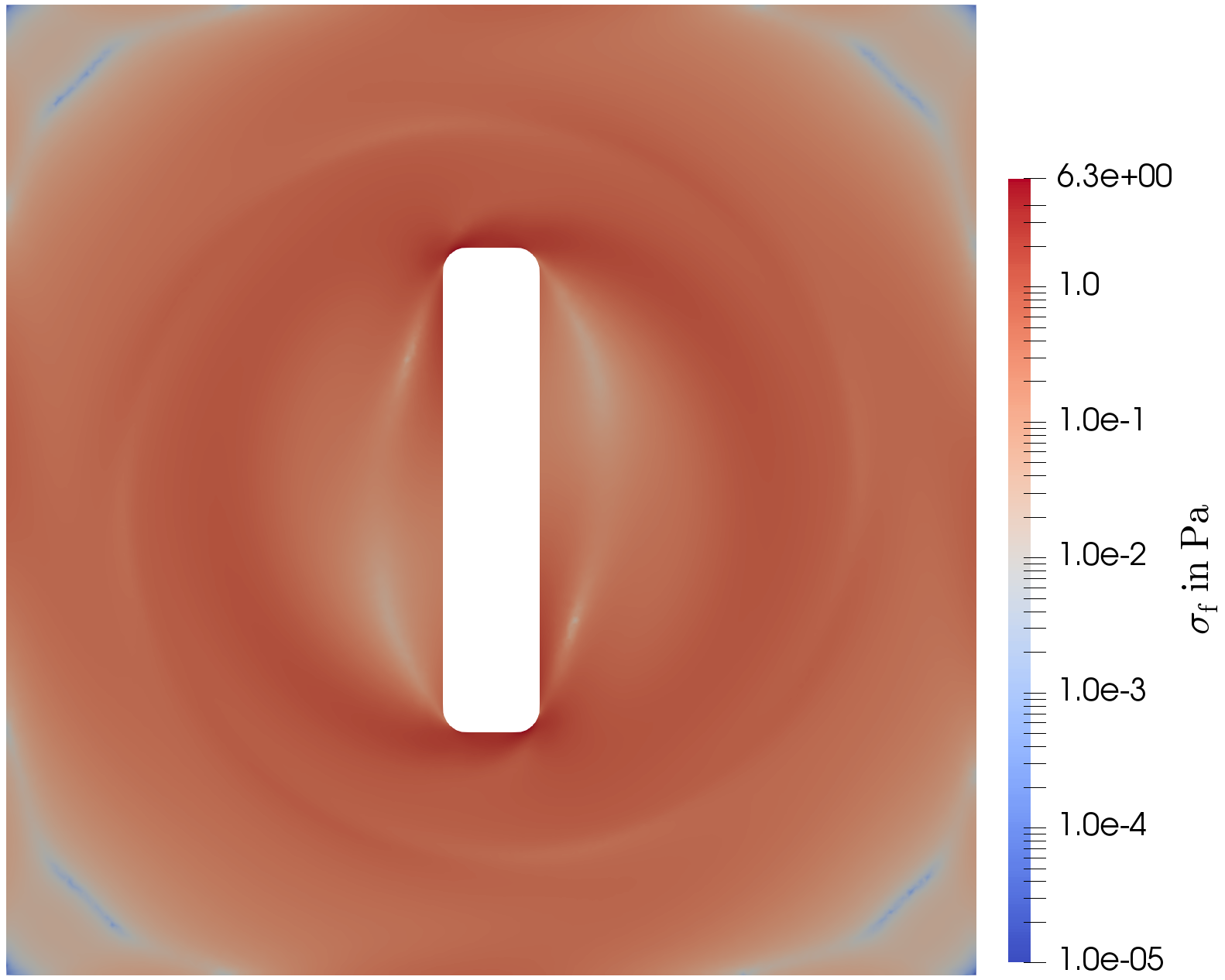}
\caption{Instantaneous shear stress $\sigma_{\mathrm{f}}$ from flow field.}
\label{fig:StirrerGf}
\end{subfigure}%
\hfill
\begin{subfigure}{.48\textwidth}
\centering
\includegraphics[width=\textwidth]{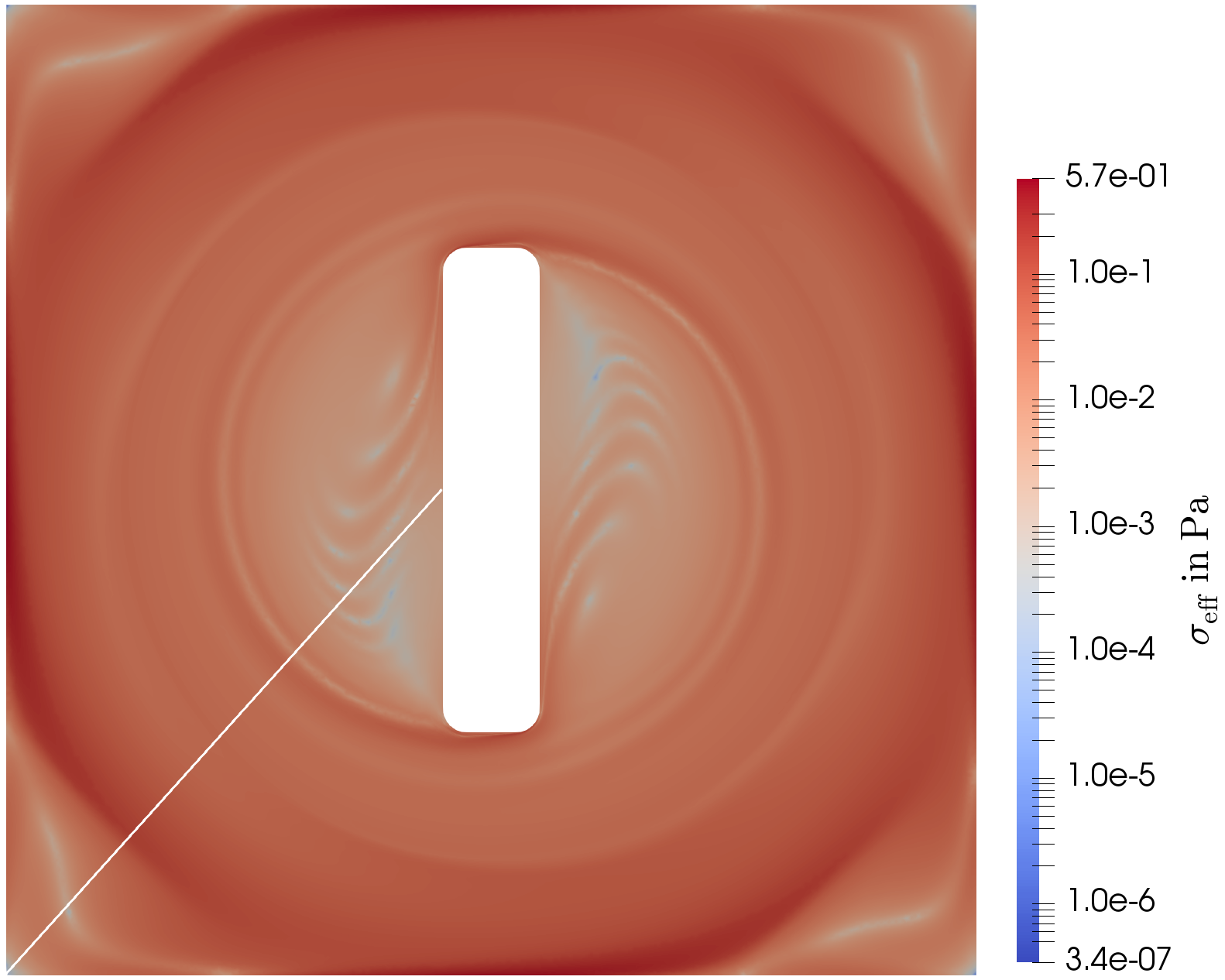}
\caption{Effective shear stress $\sigma_{\mathrm{eff}}$ from log-morph-VMS.}
\label{fig:StirrerGeff}
\end{subfigure}
\vspace{3mm}

\begin{subfigure}{\textwidth}
\centering
\includegraphics[width=.98\textwidth]{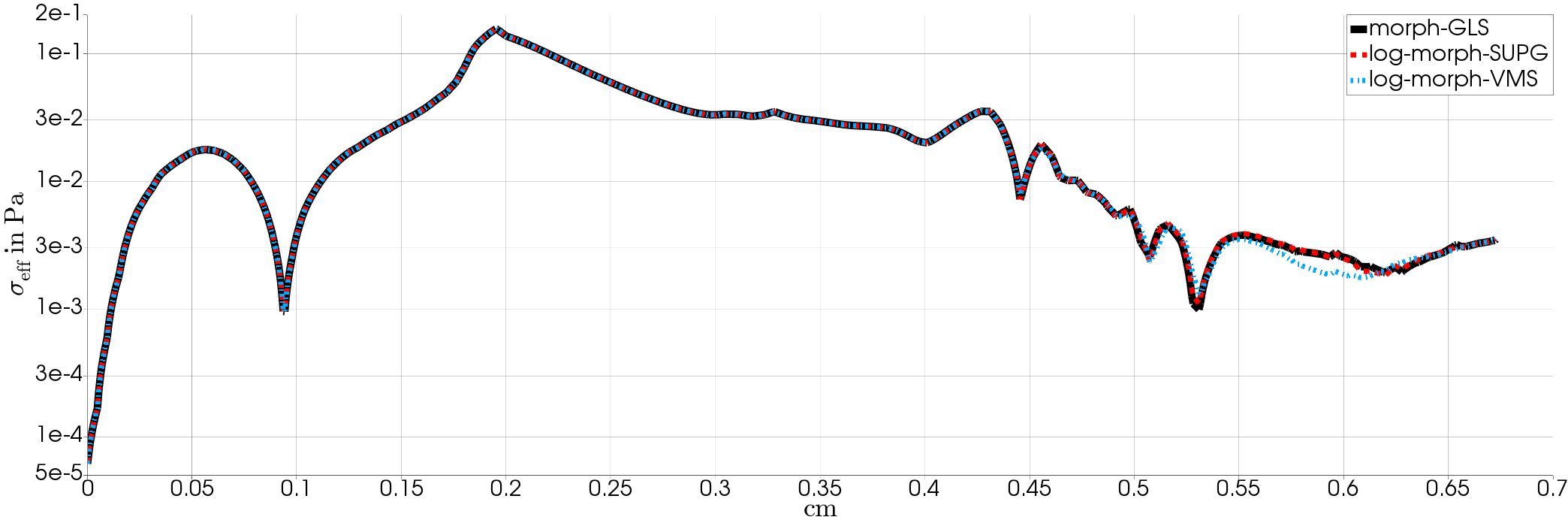}
\caption{Comparison of the effective shear stress on a line from the lower left corner to the middle of the left stirrer wall.}
\label{fig:StirrerComparison}
\end{subfigure}

\caption{Comparison of instantaneous with effective shear stress computed with three different discretizations.}
\end{figure}

A comparison of the instantaneous shear stress $\sigma_{\mathrm{f}}$ and the effective shear stress $\sigma_{\mathrm{eff}}$ is shown in Fig.~\ref{fig:StirrerGf} and Fig.~\ref{fig:StirrerGeff}.
It can be noted that the morphology equation predicts a one order of magnitude smaller shear stress acting on the RBCs compared to the instantaneous fluid stresses. As can be seen in Figure~\ref{fig:StirrerComparison}, there is hardly any difference between the morph-GLS and the log-morph-SUPG result. The effective shear stress of the log-morph-VMS simulation shows only minor differences in the inner stirrer region, where the stresses are two orders of magnitude smaller compared to the maximum value. The advantage and the enhanced numerical stability of the logarithmic transformation can be seen in Table~\ref{tbl:StirrerSolverCharacteristics}, where we observe a very significant decrease in both the GMRES and Newton-Raphson iterations needed to reach a converged solution.
\begin{table}
\centering
\begin{tabular}{l c c c c}
\hline
method         & $n_{\mathrm{GMRES}}$ & $n_{\mathrm{NR}}$ & $\varepsilon_{\det}$ & $\max|\det(\bS_i) - 1|$ \\
\hline
morph-GLS      &         2808         &         930       & $6.22\cdot10^{-5}$   & $7.67\cdot10^{-6}$\\
log-morph-SUPG &          603         &         213       & $4.23\cdot10^{-12}$  & $6.03\cdot10^{-13}$\\
log-morph-VMS  &          596         &         212       & $4.32\cdot10^{-12}$  & $4.59\cdot10^{-13}$\\
\hline
\end{tabular}
\caption{Characteristics of the simulation run for the different discretizations used for the 2D stirrer test case.}
\label{tbl:StirrerSolverCharacteristics}
\end{table}

We use the deviation of the determinant of the shape tensor from 1 over the whole domain, \ie,
\begin{equation}
\varepsilon_{\det} = \sqrt{\sum_i^{n_n} \left(\det\!\left(\bS_i\right) - 1\right)^2}
\end{equation}
as a measure for the quality of the volume conservation. Also, this measure illustrates the superiority of the log-morph formulation by a dramatic improvement of the volume conservation, as can be seen in Table~\ref{tbl:StirrerSolverCharacteristics}. For this simple test case, we cannot see any considerable differences in the log-morph-SUPG and the log-morph-VMS formulation. As a final remark, it should be mentioned that the stability of the log-morph formulation allows us to increase the time step size to values of the order of $\SI{0.1}{s}$ without any problems with convergence or the volume conservation.

\subsection{State-of-the-Art VAD}

\begin{figure}[t]
\centering
\includegraphics[width=\textwidth]{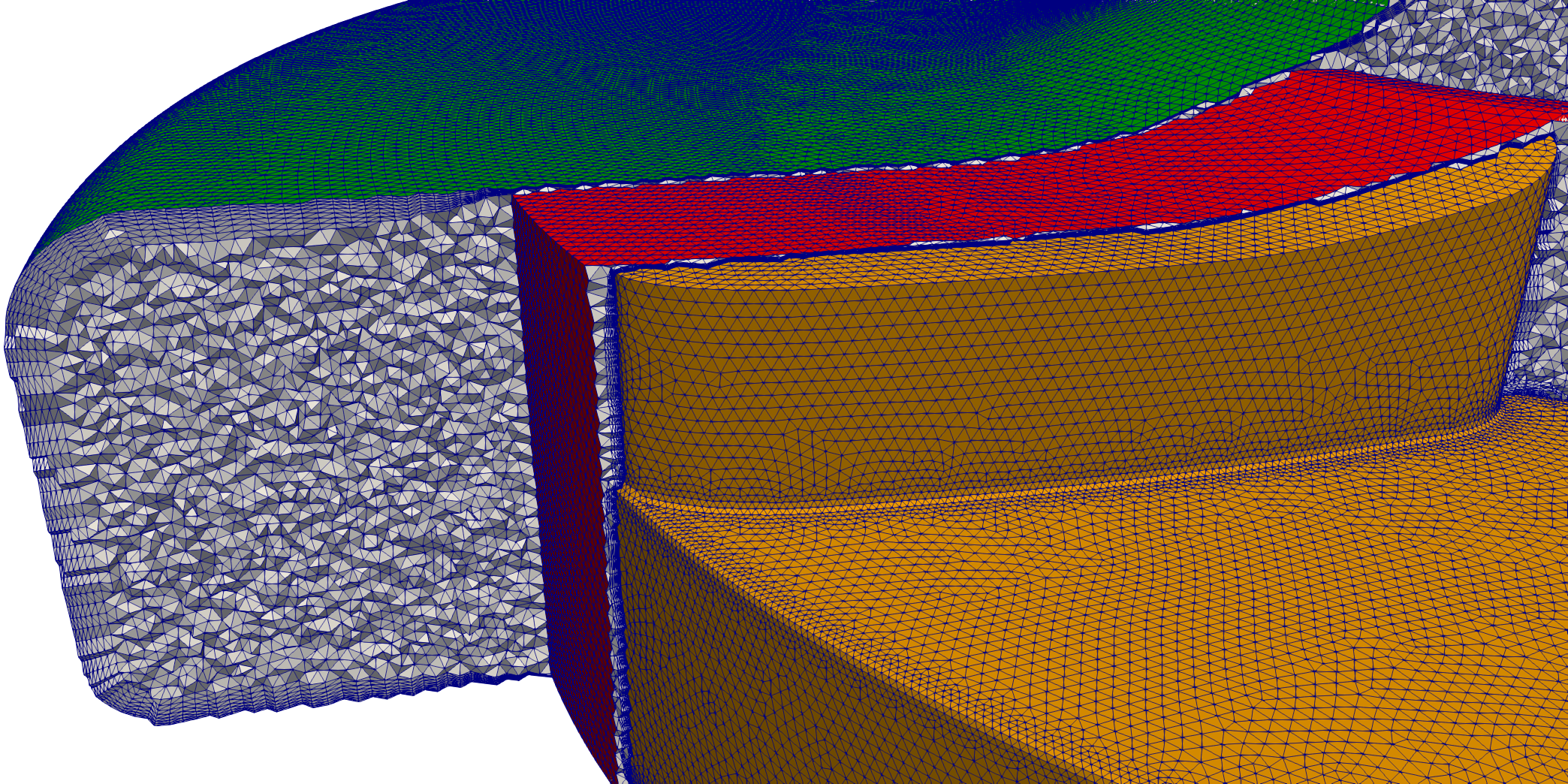}
\caption{Part of the computational mesh for the ReinVAD test case with the MRF interface (red).}
\label{fig:VADmesh}
\end{figure}

The second test case we investigate is a state-of-the-art VAD in preclinical testing developed by the ReinVAD GmbH in Aachen. The computational mesh consists of $n_n = \SI{4.74}{M}$ nodes and $n_e = \SI{27.4}{M}$ unstructured tetrahedral elements. We use a boundary layer mesh near the no-slip walls of a total thickness of $\SI{500}{\mu m}$ and seven layers with a growth rate of $\num{1.2}$, which is compressed in regions with small gap widths. We introduce an MRF interface that encloses the impeller and lies completely in the fluid volume. A part of the computational mesh together with the interface is shown in Fig.~\ref{fig:VADmesh}.

\begin{figure}[t]
\centering
\includegraphics[width=.6\textwidth]{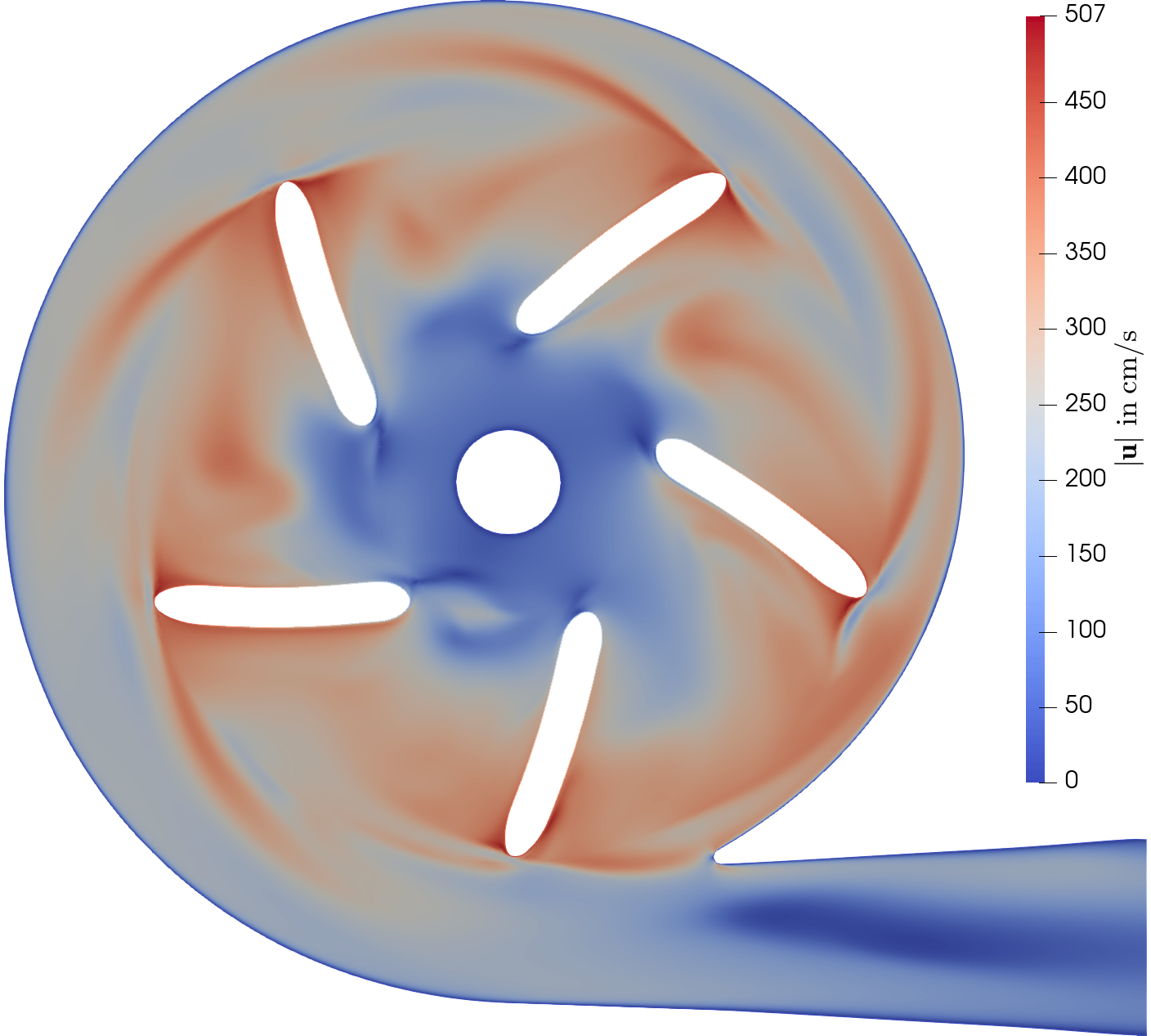}
\caption{Velocity on a slice in the middle of the impeller region.}
\label{fig:VADvelocity}
\end{figure}

We analyze the pump for an impeller angular velocity of $\SI{2400}{rpm}$ and a flow rate of $\SI{5}{L/min}$ at the inflow, using again the Newtonian blood model. We compute the steady blood flow with the commercial flow solver Altair AcuSolve using the MRF method and the SST k-$\omega$ turbulence model. The flow solution is depicted for a slice in the middle of the impeller blades in Fig.~\ref{fig:VADvelocity}.

The flow field is again used as an ambient velocity for the morphology estimation. For this complex geometry, the untransformed morphology equation, \ie, the morph-GLS method, is not able to give a converged solution. Even with the augmented Lagrangian method, a discontinuity capturing technique, and a time step size as small as $\Delta t = \SI{e-4}{s}$, the volume conservation cannot be satisfied and negative eigenvalues in the shape tensor $\bS$ occur. Nevertheless, the log-morph simulations are able to give results using a scale factor of $\alpha_\tau = 2$ for the stabilization parameter $\tau$ (eq.~\eqref{eq:tau}) and utilizing a discontinuity capturing similar to the proposition of Shakib \etal~\cite{Shakib91b}. The additional discontinuity capturing term in eq.~\eqref{eq:weak-log-morph-vms} becomes
\begin{equation}
\alpha_{\mathrm{DC}} \int_{Q_n} \nu_{\mathrm{DC}}\!\left(\res^h\right) \nabla\bPhi^h \cdot \boldsymbol{G}^{-1} \nabla\bPsi^h \diff Q,
\end{equation}
with another scale factor $\alpha_{\mathrm{DC}}$, the contravariant metric tensor $\boldsymbol{G}^{-1}$, and a numerical diffusion defined as
\begin{equation}
\nu_{\mathrm{DC}}\!\left(\res^h\right) = \sqrt{\frac{\res^h \cdot \res^h}{\nabla\bPsi^h \cdot \boldsymbol{G}^{-1} \nabla\bPsi^h}}.
\end{equation}
We choose a discontinuity capturing scale factor of $\alpha_{\mathrm{DC}} = \num{0.05}$ to ensure convergence. Furthermore, we use the MRF method and a Krylov space of $\num{50}$, an ILUT fill-in of $\num{75}$ with a threshold of $\num{e-4}$, and a time step size of $\Delta t = \SI{e-2}{s}$ and simulate for $\SI{1.5}{s}$ physical time. At the inflow, we prescribe a boundary condition for a fully developed pipe flow in axial direction (cf. Ref.~\cite{Pauli2016}).

\begin{figure}[t]
\begin{subfigure}{.48\textwidth}
\includegraphics[width=\textwidth]{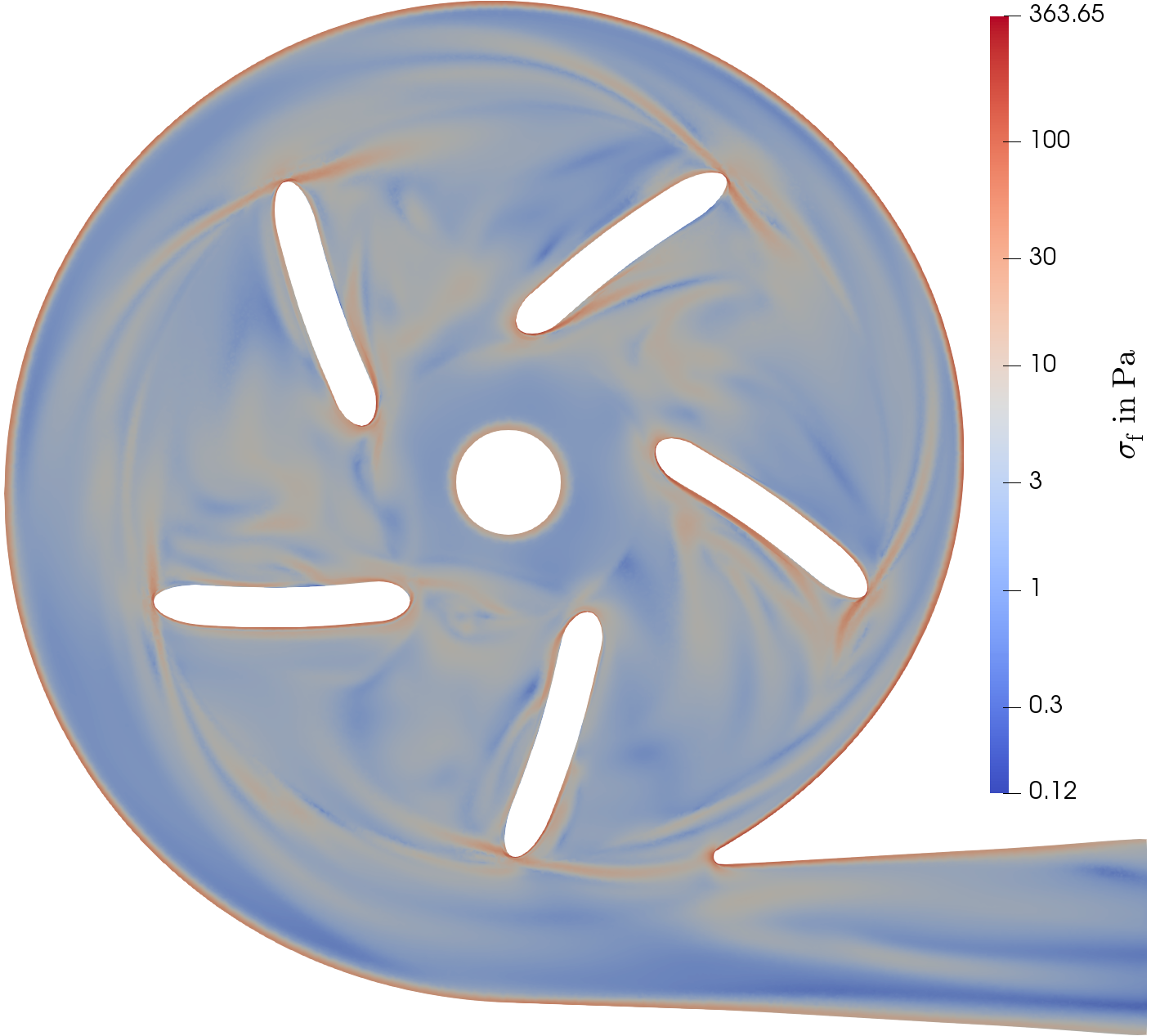}
\caption{Instantaneous shear stress $\sigma_{\mathrm{f}}$ from flow field.}
\label{fig:VAD-Gf-z}
\end{subfigure}
\hfill
\begin{subfigure}{.48\textwidth}
\includegraphics[width=\textwidth]{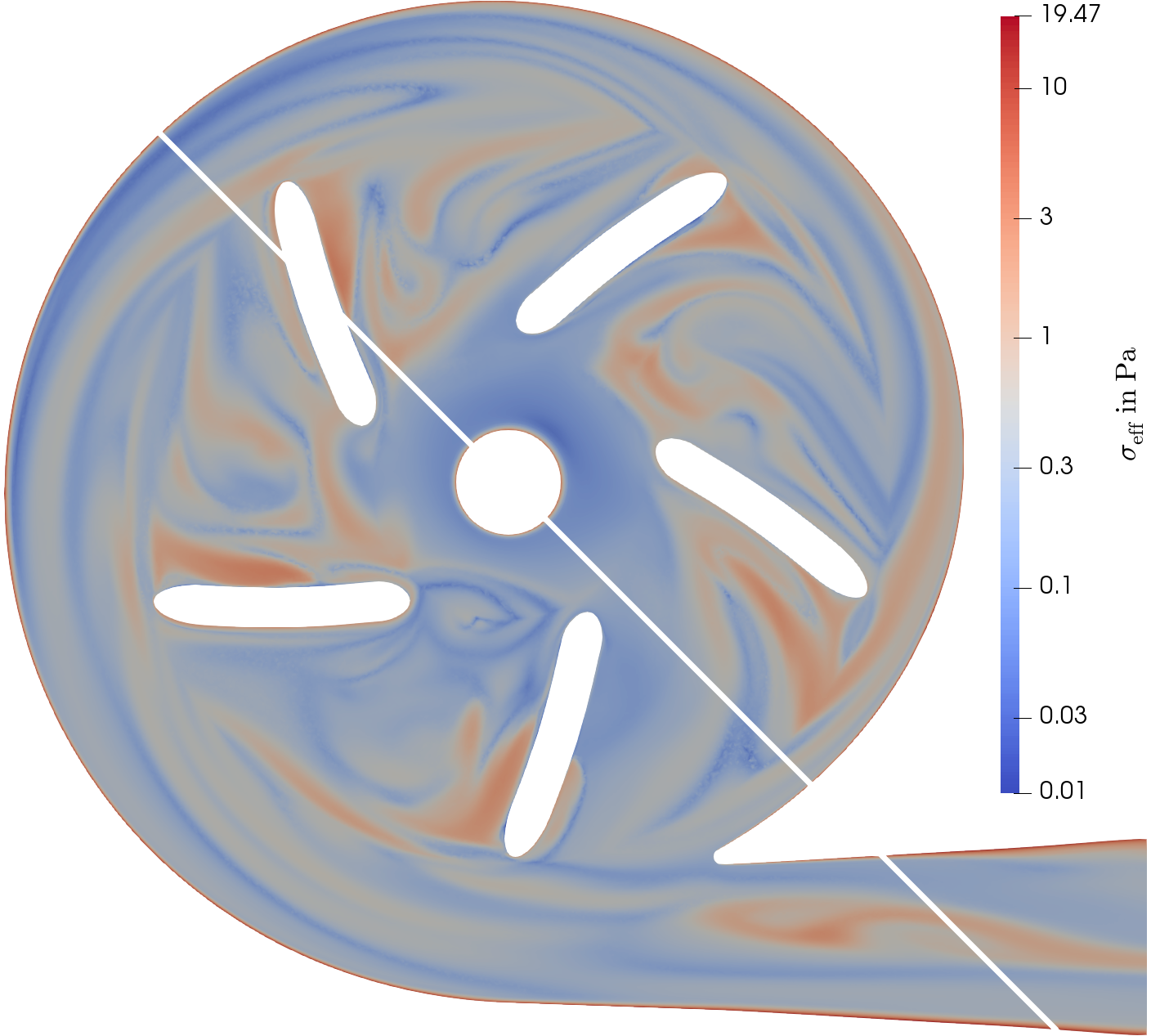}
\caption{Effective shear stress $\sigma_{\mathrm{eff}}$ from log-morph-VMS.}
\label{fig:VAD-Geff-z-VMS}
\end{subfigure}
\vspace{3mm}

\begin{subfigure}{\textwidth}
\includegraphics[width=\textwidth]{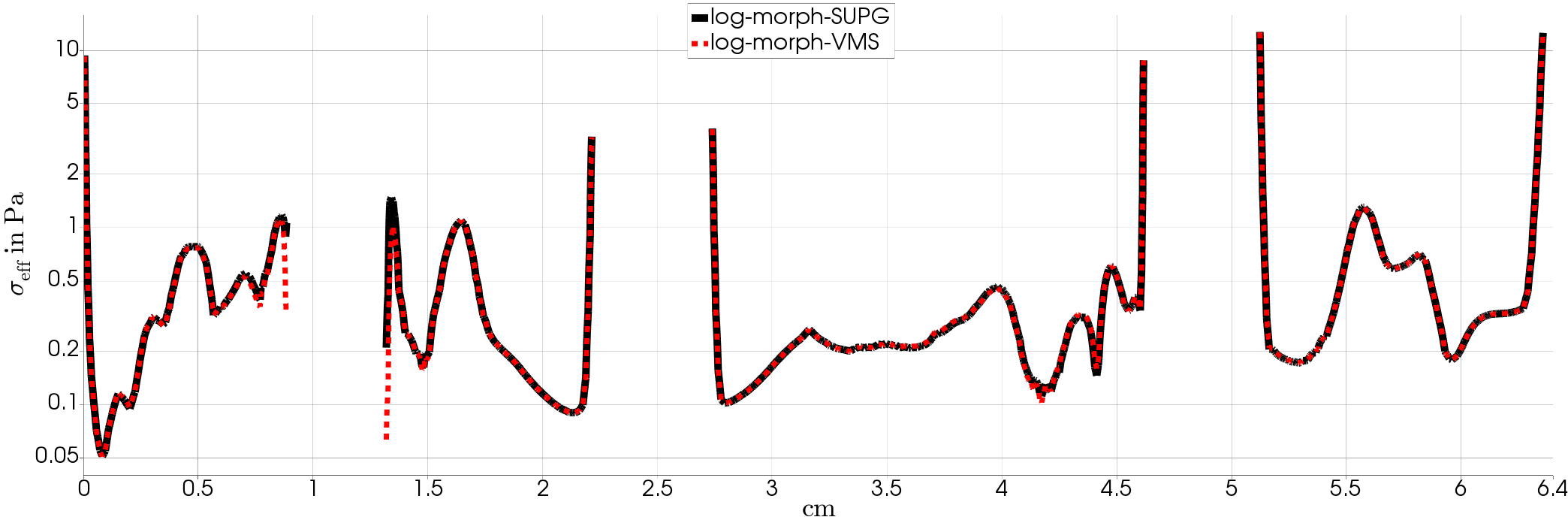}
\caption{Effective shear stress on the line shown in (b).}
\label{fig:VAD-Comparison}
\end{subfigure}
\caption{Comparison of instantaneous with effective shear stress computed with log-morph with SUPG and VMS stabilization.}
\end{figure}

Figure~\ref{fig:VAD-Gf-z} shows the instantaneous shear stress $\sigma_{\mathrm{f}}$ computed from the flow solution on a slice in the middle of the impeller region.
If compared with the effective shear stress $\sigma_{\mathrm{eff}}$ in Fig.~\ref{fig:VAD-Geff-z-VMS}, which is computed from the RBC's shape estimated by the log-morph-VMS simulation, we see a one order of magnitude lower stress. The instantaneous shear stress is especially high at the impeller blades and their tips in contrast to the effective shear stress. This behavior can be explained with the short exposure times in these regions, leading to relatively low effective stresses acting on the RBCs.

In the line plot in Fig.~\ref{fig:VAD-Comparison} the log-morph-SUPG and the log-morph-VMS method are compared. Both methods predict very similar results in most of the domain. However, close to the no-slip walls at the rotor blades (at a distance of $\SI{0.88}{cm}$ and $\SI{1.32}{cm}$) the log-morph-VMS formulation predicts significantly smaller shear stresses compared to the log-morph-SUPG method.

During our simulations, we found convergence issues for the log-morph-SUPG method. The residual fluctuated around a value of $\num{3.7e-10}$ and hence, did not reach the residual threshold set to $\num{e-10}$, always using the maximum number of Newton-Raphson iterations of 12 per time step. Another indication for convergence issues is the rather poor volume conservation shown in Table~\ref{tbl:VADSolverCharacteristics}.
\begin{table}[b]
\centering
\begin{tabular}{l c c c c}
\hline
method         & $n_{\mathrm{GMRES}}$ & $n_{\mathrm{NR}}$ & $\varepsilon_{\det}$ & $\max|\det(\bS_i) - 1|$ \\
\hline
log-morph-SUPG &          3600        &        1800       &       $0.169$        & $0.0880$ \\
log-morph-VMS  &          1180        &         590       &  $2.74\cdot10^{-8}$  & $1.47\cdot10^{-9}$\\
\hline
\end{tabular}
\caption{Characteristics of the simulation run for the different stabilizations of the log-morph method for the ReinVAD test case.}
\label{tbl:VADSolverCharacteristics}
\end{table}
In contrast, the log-morph-VMS method shows a good convergence and volume conservation using the same solver parameters.

\begin{figure}[p]
\begin{subfigure}{\textwidth}
\includegraphics[height=.3\textheight]{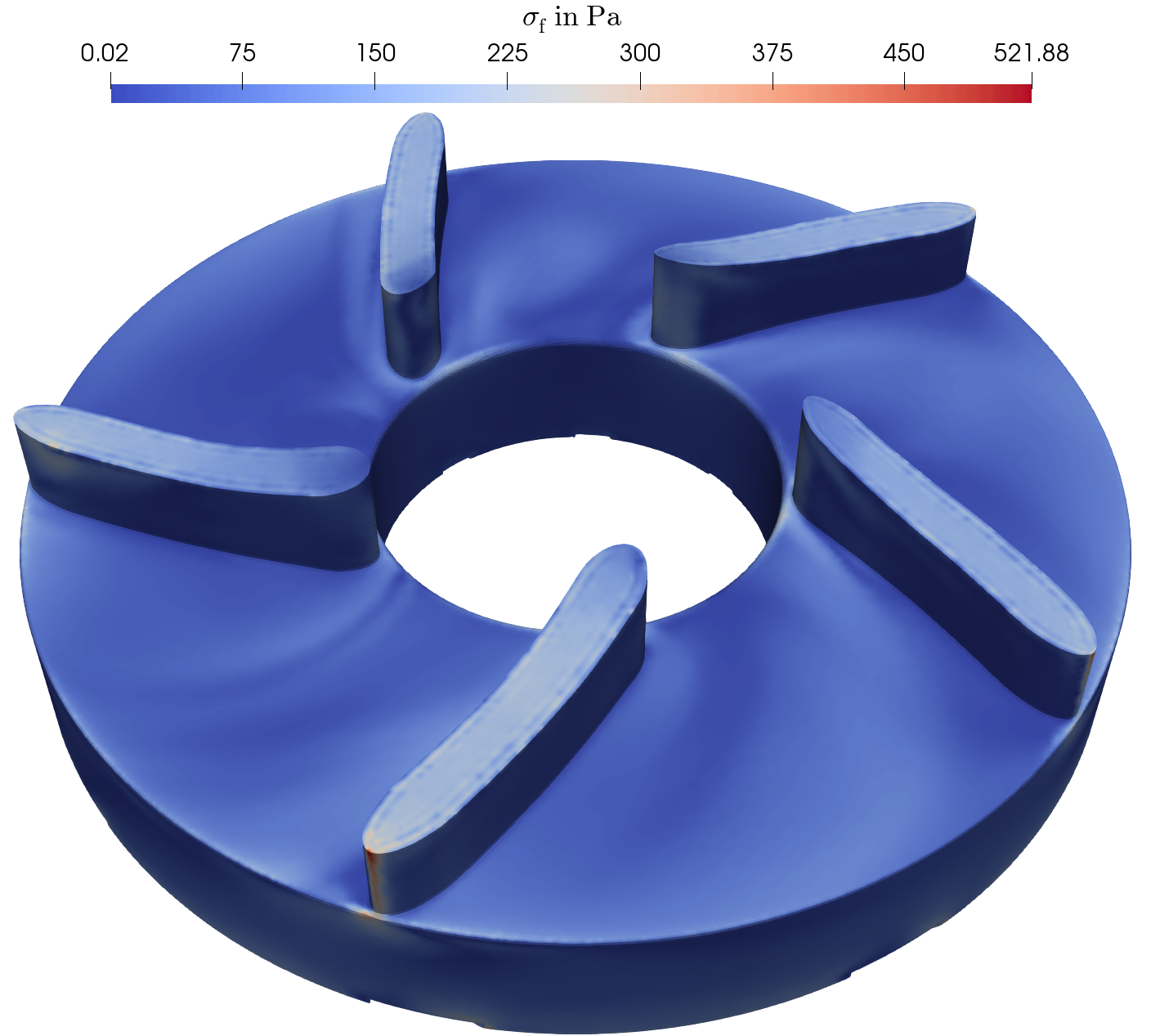}
\hfill
\includegraphics[height=.3\textheight]{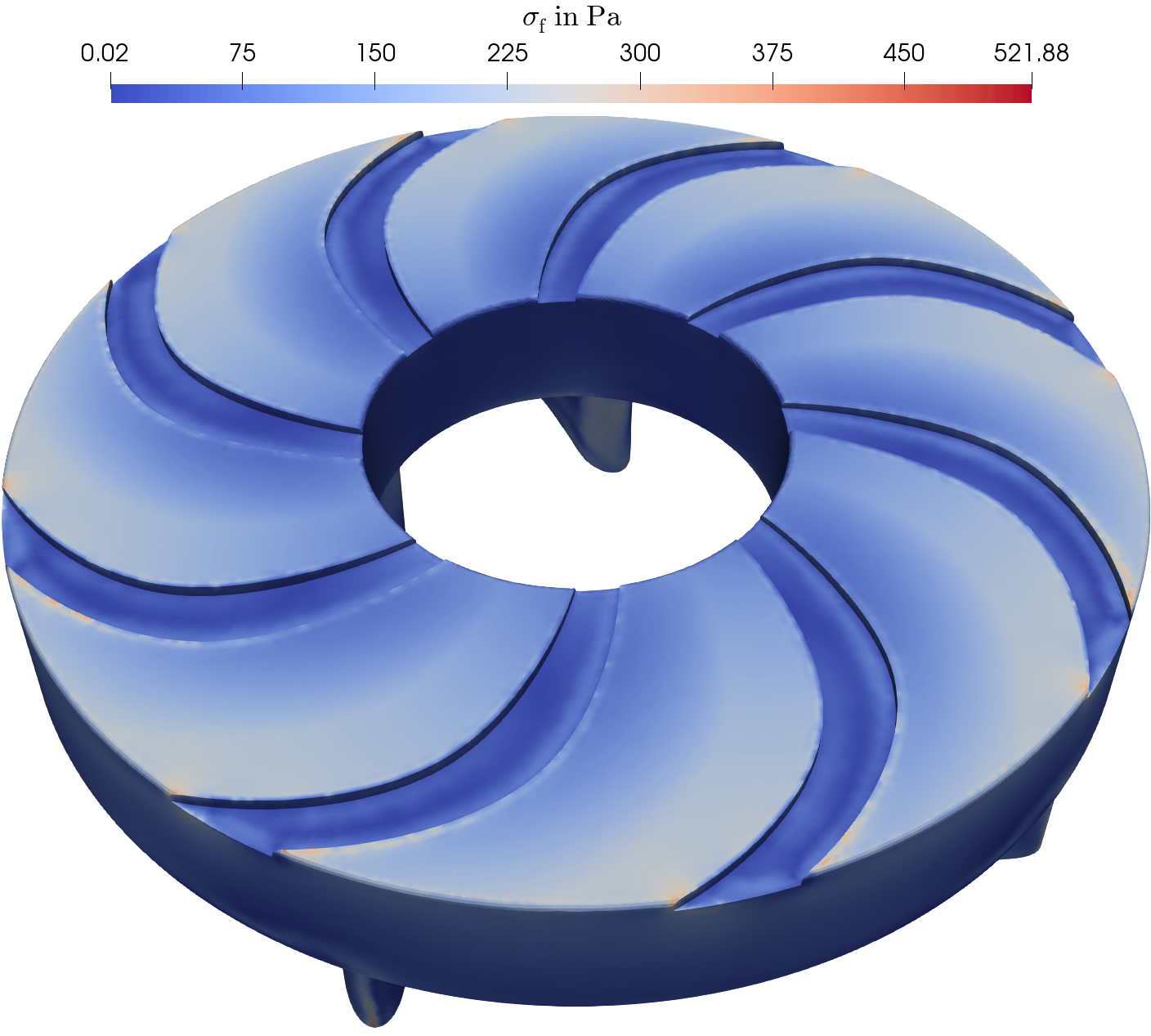}
\caption{Instantaneous shear stress $\sigma_{\mathrm{f}}$.}
\end{subfigure}
\begin{subfigure}{\textwidth}
\includegraphics[height=.3\textheight]{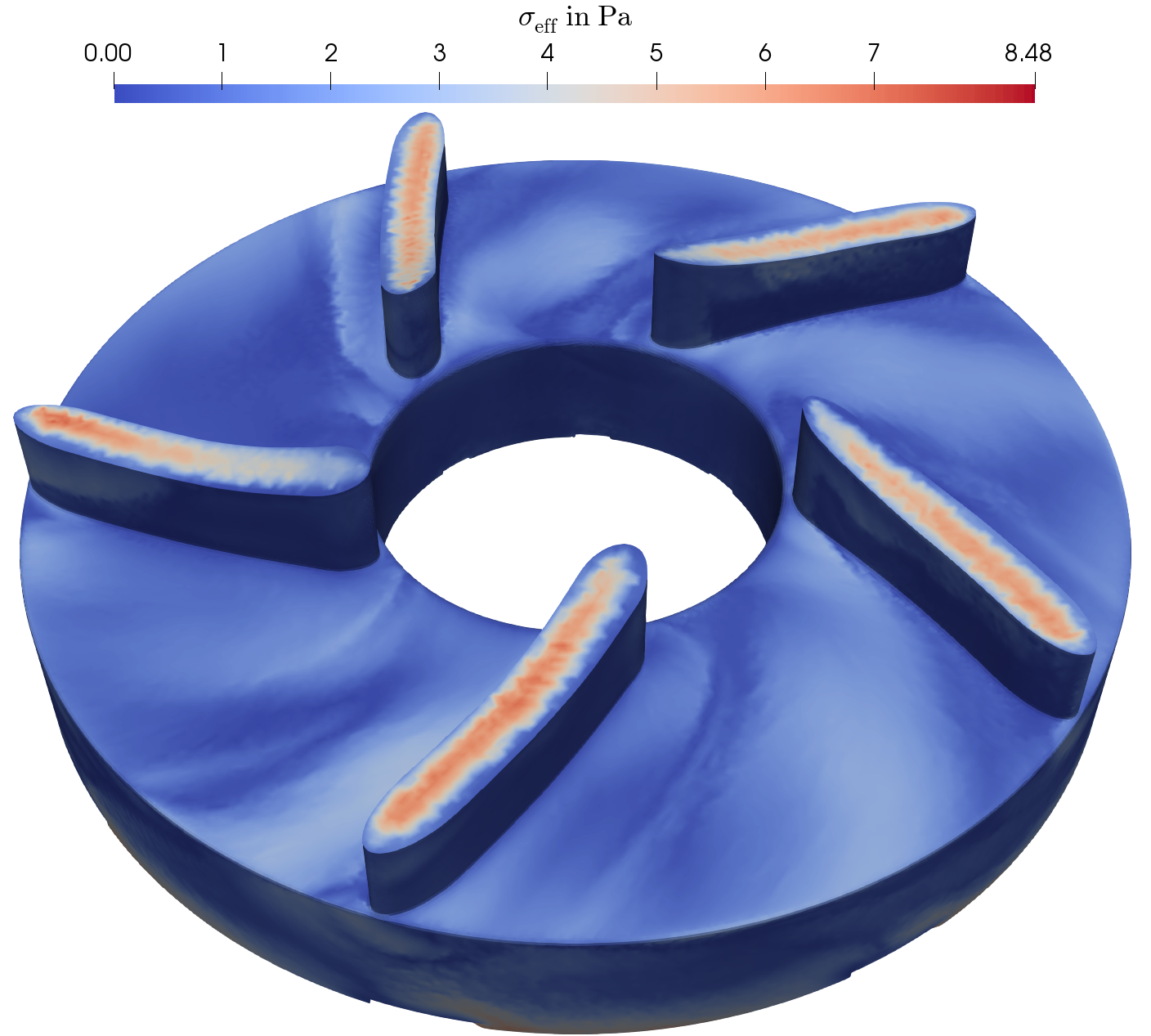}
\hfill
\includegraphics[height=.3\textheight]{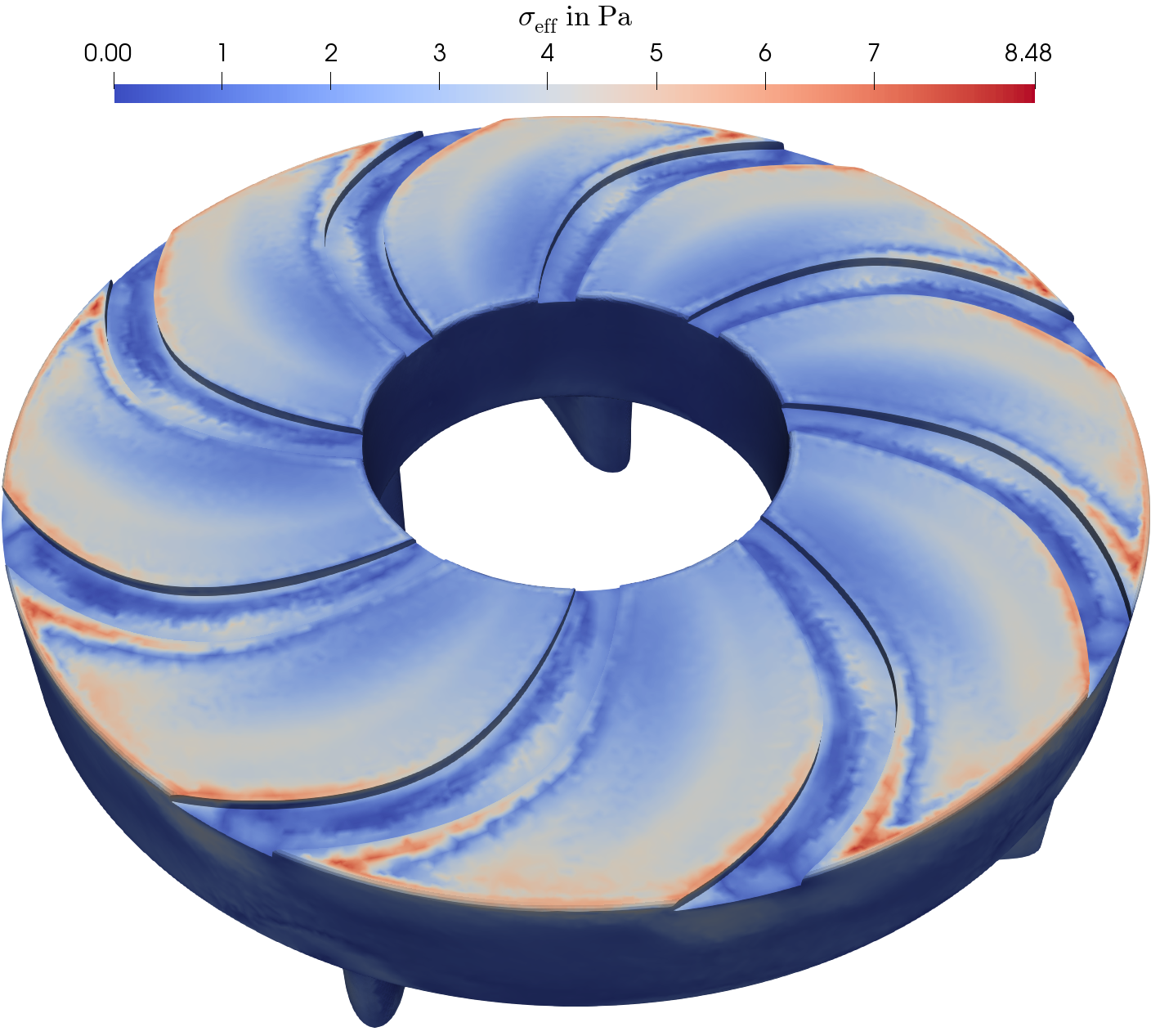}
\caption{Effective shear stress $\sigma_{\mathrm{eff}}$ from log-morph-SUPG.}
\end{subfigure}
\begin{subfigure}{\textwidth}
\includegraphics[height=.3\textheight]{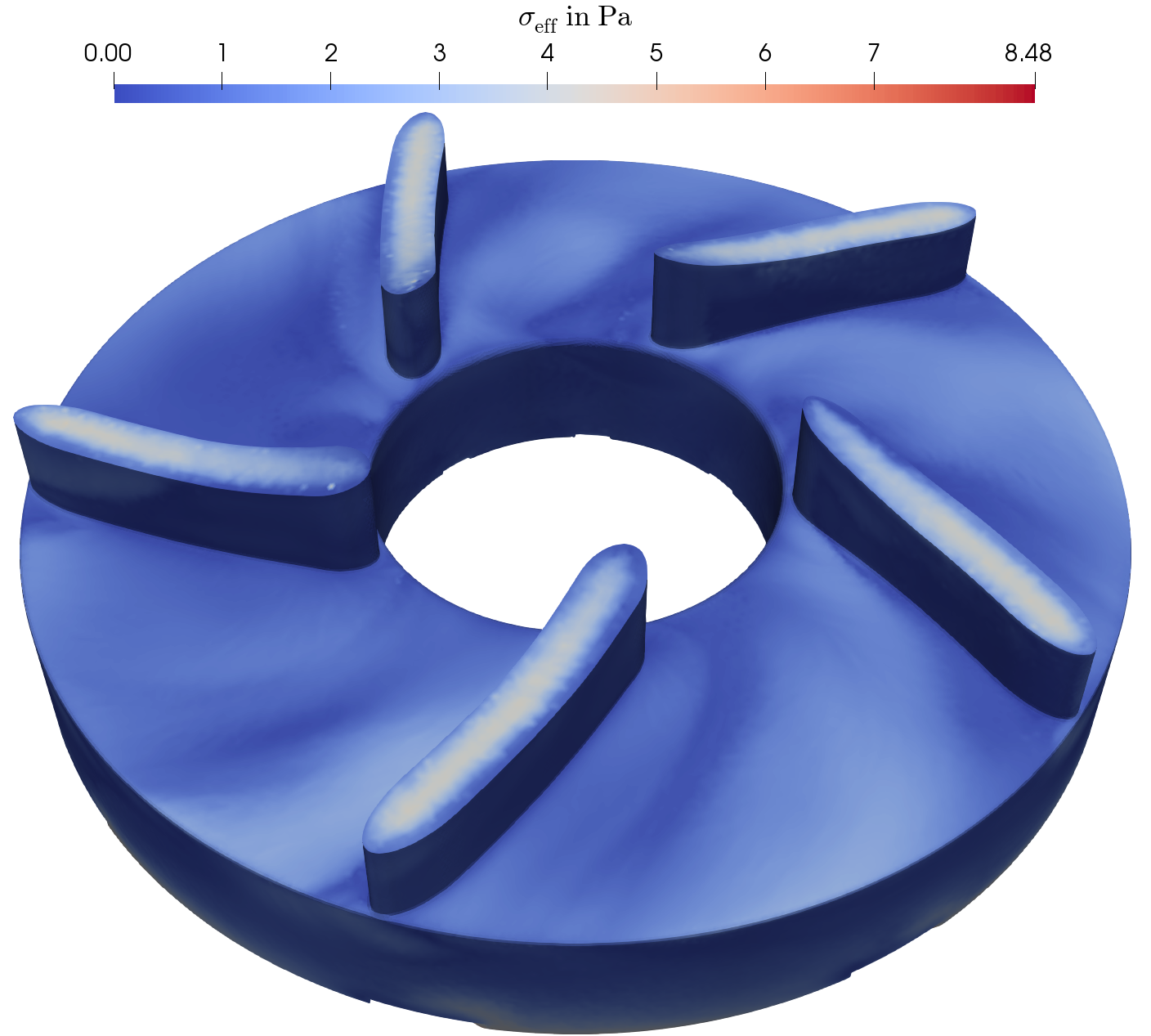}
\hfill
\includegraphics[height=.3\textheight]{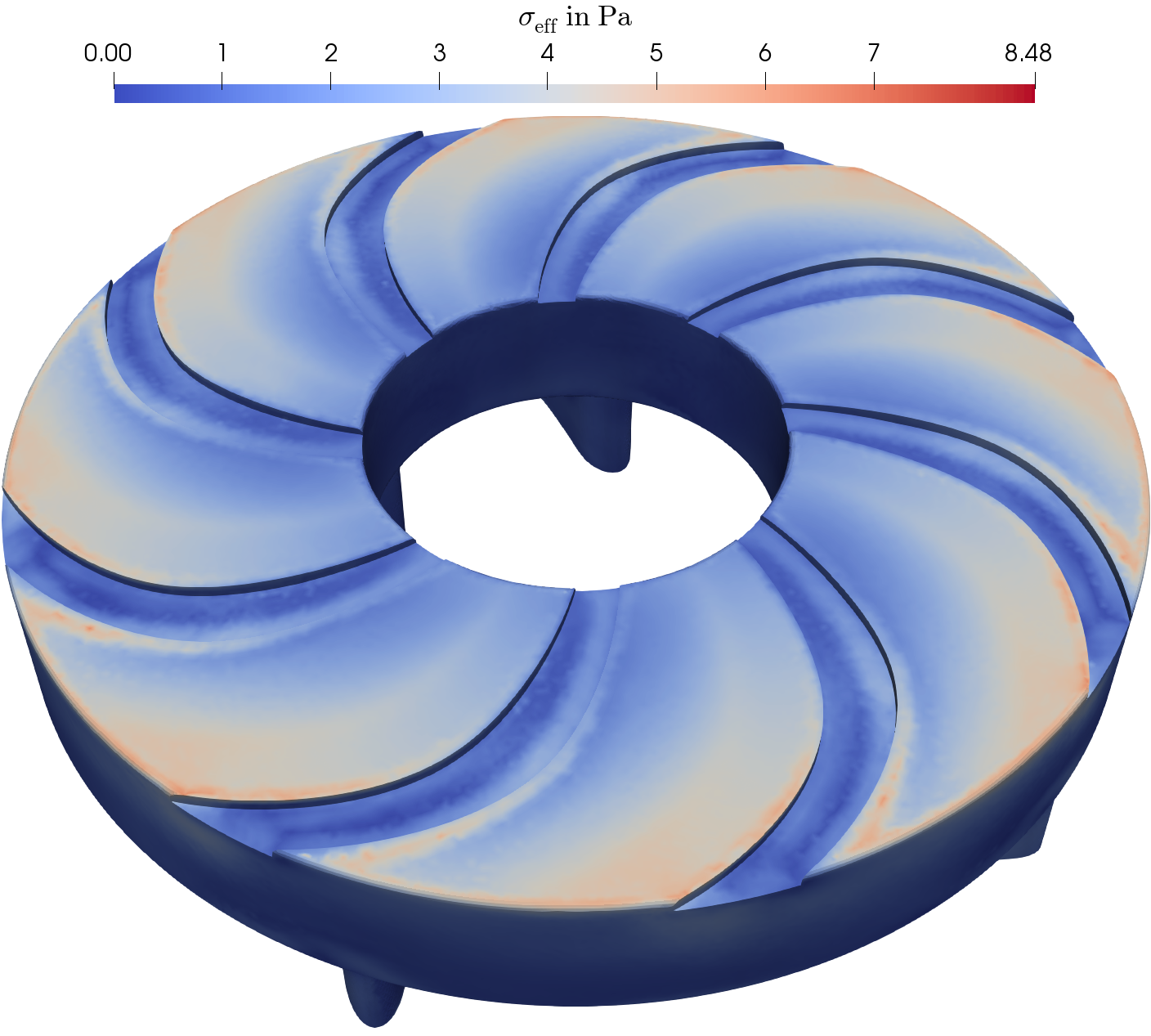}
\caption{Effective shear stress $\sigma_{\mathrm{eff}}$ from log-morph-VMS.}
\label{fig:VAD-Impeller-VMS}
\end{subfigure}
\caption{Comparison of instantaneous and effective shear stresses on the impeller top and bottom surfaces.}
\label{fig:VAD-Impeller}
\end{figure}

Another interesting comparison is the distribution of the instantaneous shear stress $\sigma_{\mathrm{f}}$ and the effective shear stress $\sigma_{\mathrm{eff}}$ for the two different discretizations on the impeller. Fig.~\ref{fig:VAD-Impeller} depicts the distribution on top and bottom of the impeller.
It can be noted that the distribution of instantaneous and effective shear stress shows similarities, such as the distribution pattern on the impeller table or in the channels of the hydrodynamic bearing. However, also clear differences are observed, especially at the impeller tips, where the peak values of $\sigma_{\mathrm{f}}$ are not present for $\sigma_{\mathrm{eff}}$. In general, the morphology equation predicts stress values about two orders of magnitude smaller than the instantaneous fluid stresses.

For the discontinuity capturing, we can only use an $\alpha_\mathrm{DC} = 0.05$ to obtain convergence for the transient log-morph simulations, although we believe that a higher scale factor $\alpha_\mathrm{DC}$ is necessary to reduce all oscillations at sharp inter-element discontinuities. However, for the converged log-morph-VMS solution, we are able to compute a restart solution for one more time step with a discontinuity capturing scale factor of $\alpha_\mathrm{DC} = \num{1.0}$, which adds a little bit more numerical diffusion and further decreases peak values. For the log-morph-SUPG simulation, this restart technique did not lead to convergence. From Fig.~\ref{fig:VAD-Impeller-VMS}, it can be observed that log-morph-VMS results in a smoother effective shear stress distribution with smaller peak values compared to the log-morph-SUPG case.  Thus, the VMS formulation helps to reduce oscillations near the impeller walls. Although we present the effective shear stress distribution for the restarted log-morph-VMS method, the log-morph-VMS result before the restart shows already fewer oscillations than the log-morph-SUPG method.

\section{Discussion}

For the simple stirrer test case, the numerical behavior of the log-morph-SUPG and the log-morph-VMS method show a comparable performance. For the complex VAD test case, though, the VMS-stabilized log-morphology formulation reveals superior convergence behavior and superior stabilization properties at the rotating impeller compared to the SUPG variant. Evidently, this behavior is due to the missing source term stabilization of the SUPG method. The rapidly decreasing velocity magnitudes in the thin boundary layer elements close to the no-slip impeller walls lead to dominating source term contributions. Hence, the VMS-stabilized formulation that by construction includes a source term stabilization is advantageous.

Although the presented VMS method leads to a significant improvement, further developments would be desirable. We expect that an enhanced framework for discontinuity capturing could further improve the convergence rate. Together with an investigation of the definition of the stabilization parameter, this could lead to a better formulation that does not need a scaling of the stabilization terms for complex geometries. Another interesting topic would be the investigation of higher order approximations for the highly nonlinear VMS terms; here, we did only consider the first Fréchet derivatives in the linearization.

\section{Conclusion}

\enlargethispage{5mm}
The aim of this paper was to present the application of the variational multiscale formalism to a tensor-based red blood cell deformation model. We used a logarithmic transformation of the shape tensor in the underlying morphology equation to enhance the numerical stability. This introduces highly nonlinear terms in the resulting log-morph equation. The VMS method is used as a general framework for stabilization of (nonlinear) partial differential equations. Its resulting terms for the log-morph equation are derived with a Fréchet derivative approach for the linearization of the nonlinear terms. To our best knowledge, this is the first application of the VMS formulation to such a highly nonlinear tensor model.

For a simple two dimensional stirrer test case, we found good agreement between the untransformed morphology equation and the log-morph equation with SUPG and VMS stabilization. The enhanced numerical stability of the log-morph equation is clearly observed for this simple test case. Furthermore, we successfully applied the log-morph-VMS method to a state-of-the-art ventricular assist device for which the untransformed morphology equation is not able to converge. Our studies showed that the VMS stabilization method leads to an improved numerical stability for complex test cases compared to an SUPG stabilization scheme.

\subsubsection*{Acknowledgments}

We gratefully acknowledge the fruitful discussions with Philipp Knechtges, without which this paper would not have been possible. We would also like to thank the ReinVAD GmbH for providing us with the investigated VAD geometry data. The authors gratefully acknowledge the computing time granted through JARA-HPC on the supercomputer JURECA at Forschungszentrum Jülich.

\bibliography{log-morph.bib}

\begin{thebibliography}{10}

\bibitem{Kirklin2015a}
J.~Kirklin, D.~Naftel, F.~Pagani, R.~Kormos, L.~Stevenson, E.~Blume, S.~Myers,
  M.~Miller, J.~Baldwin, and J.~Young, ``Seventh {INTERMACS} annual report:
  15,000 patients and counting'', {\em The Journal of Heart and Lung
  Transplantation}, {\bf 34}~(2015)~1495--1504.

\bibitem{Fraser2011a}
K.~Fraser, M.~Taskin, B.~Griffith, and Z.~Wu, ``The use of computational fluid
  dynamics in the development of ventricular assist devices'', {\em Medical
  Engineering \& Physics}, {\bf 33}~(2011)~263--280.

\bibitem{Yu2017a}
H.~Yu, S.~Engel, G.~Janiga, and D.~Th{\'e}venin, ``A review of hemolysis
  prediction models for computational fluid dynamics'', {\em Artificial
  Organs}, {\bf 41}~(2017)~603--621.

\bibitem{Heuser80}
G.~Heuser and R.~Opitz, ``A {C}ouette viscosimeter for short time shearing of
  blood'', {\em Biorheology}, {\bf 17}~(1980)~17--24.

\bibitem{Wurzinger86a}
L.~Wurzinger, R.~Opitz, and H.~Eckstein, ``Mechanical blood trauma: an
  overview'', {\em Angeiologie}, {\bf 38}~(1986)~81--97.

\bibitem{Zhang2011a}
T.~Zhang, M.~Taskin, H.~Fang, A.~Pampori, R.~Jarvik, B.~Griffith, and Z.~Wu,
  ``Study of flow-induced hemolysis using novel {C}ouette-type blood-shearing
  devices'', {\em Artificial Organs}, {\bf 35}~(2011)~1180--1186.

\bibitem{Arora2004a}
D.~Arora, M.~Behr, and M.~Pasquali, ``A tensor-based measure for estimating
  blood damage'', {\em Artificial Organs}, {\bf 28}~(2004)~1002--1015, Errata
  in \emph{Artificial Organs}, \textbf{36} (2012) 500.

\bibitem{Pauli2012b}
L.~Pauli, J.~Nam, M.~Pasquali, and M.~Behr, ``Transient stress-based and
  strain-based hemolysis estimation in a simplified blood pump'', {\em
  International Journal for Numerical Methods in Biomedical Engineering}, {\bf
  29}~(2013)~1148--1160.

\bibitem{Chen2011a}
Y.~Chen and M.~K. Sharp, ``A strain-based flow-induced hemolysis prediction
  model calibrated by in vitro erythrocyte deformation measurements'', {\em
  Artificial Organs}, {\bf 35}~(2011)~145--156.

\bibitem{Arwatz2013a}
G.~Arwatz and A.~Smits, ``A viscoelastic model of shear-induced hemolysis in
  laminar flow'', {\em Biorheology}, {\bf 50}~(2013)~45--55.

\bibitem{Ezzeldin2015}
H.~M. Ezzeldin, M.~D. de~Tullio, M.~Vanella, S.~D. Solares, and E.~Balaras, ``A
  strain-based model for mechanical hemolysis based on a coarse-grained red
  blood cell model'', {\em Annals of Biomedical Engineering}, {\bf
  43}~(2015)~1398--1409.

\bibitem{Sohrabi2017a}
S.~Sohrabi and Y.~Liu, ``A cellular model of shear-induced hemolysis'', {\em
  Artificial Organs}, {\bf 41}~(2017)~E80--E91.

\bibitem{Fattal2004a}
R.~Fattal and R.~Kupferman, ``Constitutive laws for the matrix-logarithm of the
  conformation tensor'', {\em Journal of Non-Newtonian Fluid Mechanics}, {\bf
  123}~(2004)~281--285.

\bibitem{Knechtges2015a}
P.~Knechtges, ``The fully-implicit log-conformation formulation and its
  application to three-dimensional flows'', {\em Journal of Non-Newtonian Fluid
  Mechanics}, {\bf 223}~(2015)~209--220.

\bibitem{Hughes1995a}
T.~Hughes, ``Multiscale phenomena: Green's functions, the
  {D}irichlet-to-{N}eumann formulation, subgrid scale models, bubbles and the
  origins of stabilized methods'', {\em Computer Methods in Applied Mechanics
  and Engineering}, {\bf 127}~(1995)~387--401.

\bibitem{Hughes96a}
T.~Hughes and J.~Stewart, ``A space-time formulation for multiscale
  phenomena'', {\em Journal of Computational and Applied Mathematics}, {\bf
  74}~(1996)~217--229.

\bibitem{Bazilevs2007a}
Y.~Bazilevs, V.~Calo, J.~Cottrell, T.~Hughes, A.~Reali, and G.~Scovazzi,
  ``Variational multiscale residual-based turbulence modeling for large eddy
  simulation of incompressible flows'', {\em Computer Methods in Applied
  Mechanics and Engineering}, {\bf 197}~(2007)~173--201.

\bibitem{Kwack2017a}
J.~Kwack, A.~Masud, and K.~Rajagopal, ``Stabilized mixed three-field
  formulation for a generalized incompressible {O}ldroyd-{B} model'', {\em
  International Journal for Numerical Methods in Fluids}, {\bf
  83}~(2017)~704--734.

\bibitem{Chien70a}
S.~Chien, ``Shear dependence of effective cell volume as a determinant of blood
  viscosity'', {\em Science}, {\bf 168}~(1970)~977--979.

\bibitem{Merrill66a}
E.~W. Merrill, E.~R. Gilliland, T.~S. Lee, and E.~W. Salzman, ``Blood rheology:
  Effect of fibrinogen deduced by addition'', {\em Circulation Research}, {\bf
  18}~(1966)~437--446.

\bibitem{Qin98a}
Z.~Qin, L.-G. Durand, L.~Allard, and G.~Cloutier, ``Effects of a sudden flow
  reduction on red blood cell rouleau formation and orientation using {RF}
  backscattered power'', {\em Ultrasound in Medicine \& Biology}, {\bf
  24}~(1998)~503--511.

\bibitem{Schmid-Schoenbein69a}
H.~Schmid-Sch{\"o}nbein and R.~Wells, ``Fluid drop-like transition of
  erythrocytes under shear'', {\em Science}, {\bf 165}~(1969)~288--291.

\bibitem{Vitale2014}
F.~Vitale, J.~Nam, L.~Turchetti, M.~Behr, R.~Raphael, M.~C. Annesini, and
  M.~Pasquali, ``A multiscale, biophysical model of flow-induced red blood cell
  damage'', {\em AIChE Journal}, {\bf 60}~(2014)~1509--1516.

\bibitem{Leverett72a}
L.~Leverett, J.~Hellums, C.~Alfrey, and E.~Lynch, ``Red blood cell damage by
  shear stress'', {\em Biophysical Journal}, {\bf 12}~(1972)~257--273.

\bibitem{Lanotte2016a}
L.~Lanotte, J.~Mauer, S.~Mendez, D.~Fedosov, J.-M. Fromental, V.~Claveria,
  F.~Nicoud, G.~Gompper, and M.~Abkarian, ``Red cells dynamics morphologies
  govern blood shear thinning under microcirculatory flow conditions'', {\em
  Proceedings of the National Academy of Sciences}, {\bf
  113}~(2016)~13289--13294.

\bibitem{Knechtges2014a}
P.~Knechtges, M.~Behr, and S.~Elgeti, ``Fully-implicit log-confirmation
  formulation of constitutive laws'', {\em Journal of Non-Newtonian Fluid
  Mechanics}, {\bf 214}~(2014)~78--87.

\bibitem{Al-Mohy2009}
A.~H. Al-Mohy and N.~J. Higham, ``The complex step approximation to the
  {F}r{\'e}chet derivative of a matrix function'', {\em Numerical Algorithms},
  {\bf 53}~(2009)~133.

\bibitem{Shakib91b}
S.~Shakib, T.~Hughes, and Z.~Johan, ``A new finite element formulation for
  computational fluid dynamics: {X}.~the compressible {Euler} and
  {Navier}-{Stokes} equations'', {\em Computer Methods in Applied Mechanics and
  Engineering}, {\bf 89}~(1991)~141--219.

\bibitem{Pauli2016b}
L.~Pauli and M.~Behr, ``On stabilized space-time {FEM} for anisotropic meshes:
  Incompressible {Navier}-{Stokes} equations and applications to blood flow in
  medical devices'', {\em International Journal for Numerical Methods in
  Fluids}, {\bf 85}~(2017)~189--209.

\bibitem{Hughes98a}
T.~Hughes, G.~Feij\'{o}o, L.~Mazzei, and J.-B. Quincy, ``The variational
  multiscale method---a pradigm for computational mechanics'', {\em Computer
  Methods in Applied Mechanics and Engineering}, {\bf 166}~(1998)~3--24.

\bibitem{Pauli2016}
L.~Pauli, {\em Stabilized Finite Element Methods for Computational Design of
  Blood-Handling Devices}, Ph.D. thesis, RWTH Aachen University, 2016.

\bibitem{jureca}
{J\"{u}lich Supercomputing Centre}, ``{JURECA: Modular supercomputer at
  J\"{u}lich Supercomputing Centre}'', {\em Journal of large-scale research
  facilities}, {\bf 4}~(2018).

\bibitem{Pauli2014a}
L.~Pauli, J.~Both, and M.~Behr, ``Stabilized finite element method for flows
  with multiple reference frames'', {\em International Journal for Numerical
  Methods in Fluids}, {\bf 78}~(2015)~657--669.

\end{thebibliography}
\bibliographystyle{fine}

\appendix

\section{Derivation of the Volume Conservation Term\label{sec:VolumeConservation}}

The volume conservation for the morphology equation can be derived
using the determinant of $\bS$, \ie, $\det\!\left(\bS\right)$, which is proportional to the volume
of the described ellipsoid. This translates to
\begin{equation}
\det\!\left(\e^{\bPsi}\right)=\e^{\tr\left(\bPsi\right)},
\end{equation}
for the logarithmic shape tensor $\bPsi$. To preserve the volume
means mathematically that
\begin{equation}
\det\!\left(\bS\right)=\text{const, or }\tr\!\left(\bPsi\right)=\text{const, or }\frac{\diff}{\diff t}\tr\!\left(\bPsi\right)=0.
\end{equation}
This can be used in eq.~(\ref{eq:Res-log-morph}) and one can show
for an incompressible fluid $\left(\tr\!\left(\bE\right)=0\right)$
that
\begin{equation}
\tr\!\left(\boldsymbol{F}\!\left(\bPsi,\bE\right)\right) = \sum_{i,j=1}^{d} f\!\left(\lambda_i - \lambda_j\right) \tr\!\left(\proj_i \bE \proj_j\right) = \sum_{i=1}^{d} 2 \boldsymbol{e}_i^T \bE \boldsymbol{e}_i = 2 \tr\!\left(\bE\right) = 0,
\end{equation}
and that
\begin{equation}
\tr\!\left(\bW\bPsi - \bPsi\bW\right) = \tr\!\left(\bW\bPsi\right) - \tr\!\left(\bPsi\bW\right) = \tr\!\left(\bW\bPsi\right) - \tr\!\left(\bW\bPsi\right) = 0,
\end{equation}
which leads to
\begin{equation}
\frac{\diff}{\diff t} \tr\!\left(\bPsi\right) = \alpha_1 \left(\tr\!\left(\boldsymbol{1}\right) - g\!\left(\bPsi\right) \tr\!\left(\e^{-\bPsi}\right)\right) \overset{!}{=} 0,
\end{equation}
and finally results in
\begin{equation}
g\!\left(\bPsi\right) = \frac{d}{\tr\!\left(\e^{-\bPsi}\right)},
\end{equation}
with the number of space dimensions $d$.

\section{Numerical Treatment of the VMS Prefactors}\label{sec:Prefactors}

The second prefactor in the directional derivative of the $\boldsymbol{L}_{\alpha_2}$-term in eq.~\eqref{eq:DerivativeLa2} is numerically not accurate in the vicinity of small denominators. It can be rewritten with difference quotients using $x = \lambda_i - \lambda_l$, $y = \lambda_j - \lambda_l$, and $z = \lambda_k - \lambda_l$ as
\begin{align}
&\frac{f\!\left(x\right)}{\left(x - y\right) \left(x - z\right)} + \frac{f\!\left(y\right)}{\left(y - x\right) \left(y - z\right)} + \frac{f\!\left(z\right)}{\left(z - x\right)\left(z - y\right)}\nonumber \\
=& \frac{1}{3} \left[\frac{\coth\!\left(\frac{x}{2}\right) - \frac{2}{x} - \coth\!\left(\frac{y}{2}\right) + \frac{2}{y}}{x - y}
+ \frac{\coth\!\left(\frac{y}{2}\right) - \frac{2}{y} - \coth\!\left(\frac{z}{2}\right) + \frac{2}{z}}{y - z}
+ \frac{\coth\!\left(\frac{z}{2}\right) - \frac{2}{z} - \coth\!\left(\frac{x}{2}\right) + \frac{2}{x}}{z - x}\right. \nonumber \\
&\qquad +\frac{x}{y - z} \left(\frac{\coth\!\left(\frac{y}{2}\right) - \frac{2}{y} - \coth\!\left(\frac{x}{2}\right) + \frac{2}{x}}{y - x} - \frac{\coth\!\left(\frac{z}{2}\right) - \frac{2}{z} - \coth\!\left(\frac{x}{2}\right) + \frac{2}{x}}{z - x}\right) \nonumber \\
&\qquad +\frac{y}{z - x} \left(\frac{\coth\!\left(\frac{z}{2}\right) - \frac{2}{z} - \coth\!\left(\frac{y}{2}\right) + \frac{2}{y}}{z - y} - \frac{\coth\!\left(\frac{x}{2}\right) - \frac{2}{x} - \coth\!\left(\frac{y}{2}\right) + \frac{2}{y}}{x - y}\right) \nonumber \\
&\qquad +\left.\frac{z}{x - y} \left(\frac{\coth\!\left(\frac{x}{2}\right) - \frac{2}{x} - \coth\!\left(\frac{z}{2}\right) + \frac{2}{z}}{x - z} - \frac{\coth\!\left(\frac{y}{2}\right) - \frac{2}{y} - \coth\!\left(\frac{z}{2}\right) + \frac{2}{z}}{y - z}\right)\right],
\end{align}
where the difference quotients are approximated by the Taylor series \eqref{eq:TaylorSmallDenominator} for values $|x-y| < 10^{-2}$ and the arising derivatives are replaced by their Taylor series up to fourth order for $\left|\nicefrac{(x+y)}{2}\right| < 10^{-1}$.

The third prefactor can be approximated by applying the Taylor expansion twice when both denominators are small, yielding
\begin{gather}
\frac{1}{\lambda_i - \lambda_j} \left[\frac{f\!\left(\lambda_k - \lambda_i\right) - f\!\left(\lambda_l - \lambda_i\right)}{\lambda_k - \lambda_l} - \frac{f\!\left(\lambda_k - \lambda_j\right) - f\!\left(\lambda_l - \lambda_j\right)}{\lambda_k - \lambda_l}\right] = -f^{(2)}\!\left(\frac{\lambda_k + \lambda_l - \lambda_i - \lambda_j}{2}\right) \nonumber \\
-f^{(4)}\!\left(\frac{\lambda_k + \lambda_l - \lambda_i - \lambda_j}{2}\right) \left[\frac{\left(\lambda_i - \lambda_j\right)^{2}}{24} + \frac{\left(\lambda_k - \lambda_l\right)^{2}}{24}\right] + \mathcal{O}\!\left(\left(\lambda_i - \lambda_j\right)^{4}\right) + \mathcal{O}\!\left(\left(\lambda_k - \lambda_l\right)^{4}\right),
\end{gather}
which is again used when $\left|\lambda_i - \lambda_j\right| < 10^{-2}$ and where the derivatives are approximated for $\left|\nicefrac{\left(\lambda_k + \lambda_l - \lambda_i - \lambda_j\right)}{2}\right| < 10^{-1}$.

\end{document}